\documentclass[oneside,leqno,11pt]{article}
\usepackage{revstat}
\begin{document}
\title{M-ESTIMATION METHOD BASED ASYMMETRIC OBJECTIVE FUNCTION
\footnote{}}
\renewcommand{\titleheading}
             {Using the ``revstat.sty" Package}  
\author{\authoraddress{Mehmet Niyazi \c{C}ankaya}
                      {Faculty of Arts and Science, Department of Statistics,
                       U\c{s}ak University, U\c{s}ak, \\
                      Turkey
                       \ (mehmet.cankaya@usak.edu.tr)}
\\
        \authoraddress{Olcay Arslan}
                      {Faculty  Science, Department of Statistics,
                       Ankara  University, Ankara,\\
                       Turkey
                       \ (oarslan@ankara.edu.tr)}
}
\renewcommand{\authorheading}
             {Mehmet Niyazi \c{C}ankaya \ and \ Olcay Arslan}  

\maketitle

\begin{abstract}
The asymmetric $\rho$ function is proposed as an alternative to Huber $\rho$ function to model skewness and obtain robust estimators for the location, scale and skewness parameters. The robustness and asymptotic properties of the asymmetric M-estimators are explored.  A simulation study and real data examples are given to illustrate the performance of proposed asymmetric M-estimation method over the symmetric M-estimation method. It is observed from the simulation results that the asymmetric M-estimators perform better than Huber M-estimators when the data have skewness. The application on regression is also considered.
\end{abstract}



\mainpaper  


\section{Introduction} \label{intro}
The robust estimation method for the location parameter was proposed by \cite{Hub64}. In the robustness, there are different influence functions have been considered. \cite{Andetalpsicollections72,Hametal86} gave the influence functions to estimate the location and scale parameters robustly. It can be observed that these functions are symmetric. Then, it is not possible to model the asymmetry in data set. In our proposal, we will consider the asymmetric form of influence functions. The asymmetric form of Huber M-estimation will be proposed. The benefit of our asymmetric Huber M-estimation is that the location, scale and skewness parameters can be estimated when the asymmetrically data set does not come from a distribution. In other words, it is well known that the distribution assumption on data set can be a restrictive. In such a case, the location, scale and skewness parameters can be estimated by means of the function we proposed. These approaches for estimating the parameters are in the robust methods.

To get the asymmetric objective function denoted by $\rho_{ESH}$, the distributions proposed by \cite{Elsal05,Elsal08,Elsaletal05,Canetal15} will be used. A family of these distributions are proposed by \cite{Canetal15}. The some special values of parameters in distribution proposed by \cite{Canetal15} give the epsilon-skew normal (ESN) and epsilon skew Laplace distributions (ESL). The details of how one can get ESN and ESL are given by  \cite{Canetal15}. We will use these distributions to get the asymmetric objective function. In this context, we will give the following approach:

Let $f$ be a probability density function. $\rho=-log(f)$ is known to be objective function in the robustness. The normalizing constant in $f$ can be removed. Let
\begin{itemize}
  \item $\rho_{ESN}(u) = \frac{u^2}{2(1-sign(u)\varepsilon)^2}$ \\
  \\
\noindent and
\\
  \item $\rho_{ESL}(u) = \frac{|u|}{2^{1/2}(1-sign(u)\varepsilon)}$
\end{itemize}

\noindent be an objective functions of ESN and ESL distributions.

Huber's $\rho$ function is given by the following form:

\begin{equation}\label{Huberrho}
 \rho(u)=
    \left\{
  \begin{array}{ll}
    u^2~ &,  |u| \leq k; \\
    2k|u|-k^2 &, |u| > k,
  \end{array}
\right.
\end{equation}
\noindent and $\rho^{'}(u)=2\psi(u)$, that is, the function $\psi$ is a derivative of $\rho$.
\begin{equation}\label{Huberpsi}
 \psi(u)=
    \left\{
  \begin{array}{ll}
    u &, |u| \leq k; \\
   sign(u)k &, |u| > k.
  \end{array}
\right.
\end{equation}

\section{M-estimation based on Asymmetric Objective Function}

We used the Huber's $\rho$ function to propose asymmetric Huber M-estimation. The parameter $k$ in Huber is the tuning parameter to get the robust estimators. In our case, asymmetric Huber M-estimation has $c_1$ and $c_2$ that are the tuning parameters due to fact that we will propose the asymmetric form of Huber's $\rho$ function.

As it is given, the asymmetric Huber M-estimation can be considered as follow:
\begin{equation}\label{boundwithcepsilonrhokabaca}
 \rho^{*}(u)= \left\{
                \begin{array}{ll}
                  \frac{u}{2^{1/2}(1+\varepsilon)}, & (-\infty,c_1); \\
                  \frac{u^2}{2(1+\varepsilon)^2}, & [c_1,0); \\
                  \frac{u^2}{2(1-\varepsilon)^2}, & [0,c_2]; \\
                  \frac{u}{2^{1/2}(1-\varepsilon)}, & (c_2,\infty).
                \end{array}
              \right.
\end{equation}
\noindent The proposed asymmetric $\rho^{*}$ function is not continuous at  points $c_1$ and $c_2$.  After the required regularization on function $\rho^{*}$ is done, the following new $\rho^{**}$ function can be proposed as follow:
\begin{equation}\label{boundwithcepsilonrho}
 \rho^{**}(u)= \left\{
                \begin{array}{ll}
                  \frac{c_1u}{(1+\varepsilon)^2}-\frac{c_1^2}{2(1+\varepsilon)^2}, & (-\infty,c_1); \\
                  \frac{u^2}{2(1+\varepsilon)^2}, & [c_1,0); \\
                  \frac{u^2}{2(1-\varepsilon)^2}, & [0,c_2]; \\
                  \frac{c_2u}{(1-\varepsilon)^2}-\frac{c_2^2}{2(1-\varepsilon)^2}, &(c_2,\infty).
                \end{array}
              \right.
\end{equation}

\noindent where $c_1$ and $c_2$ are the tuning parameters and the continuity of $\rho^{*}$ is guaranteed via these parameters. Here, $\varepsilon$ is a skewness parameter to model the asymmetry.

\begin{definition}
The function $\rho^{**}$ in equation \eqref{boundwithcepsilonrho} is defined to be asymmetric $\rho_{ESH}$ function.
\end{definition}

$\rho_{ESH}$ is used to show the asymmetric $\rho^{**}$ function we proposed. When $c_1=c_2$ and $\varepsilon=0$, $\rho_H$ in equation \eqref{Huberrho} can be obtained.  

The function $\psi$ is a derivative of function $\rho_{ESH}$. It can be given in the following form:
\begin{equation}\label{boundwitckepsilonpsi}
 \psi(u)= \left\{
                \begin{array}{ll}
                  \frac{c_1}{(1+\varepsilon)^2}, & (-\infty,c_1); \\
                  \frac{u}{(1+\varepsilon)^2}, & [c_1,0); \\
                  \frac{u}{(1-\varepsilon)^2}, & [0,c_2]; \\
                  \frac{c_2}{(1-\varepsilon)^2}, & (c_2,\infty).
                \end{array}
              \right.
\end{equation}

The estimators of parameters $\theta$, $\sigma$ ve $\varepsilon$ can be obtained by means of asymmetric objective function given in equation \eqref{boundwithcepsilonrho}. The functions $\psi_{ESH}$, $\psi_{ESN}$ ve $\psi_{ESL}$ can be obtained from the objective functions $\rho_{ESH}$, $\rho_{ESN}$ ve $\rho_{ESL}$. When $\varepsilon=0$, the influence functions ($\psi$), the symmetric influence functions  $\psi_H,\psi_N$ and $\psi_L$ are obtained.

\subsection{M-estimators generated by asymmetric M-objective function}\label{asimetrikrho}

Suppose that the random variables $X_1,X_2,...,X_n$ are distributed as a probability density function $f$. The parameters $\theta$, $\sigma$ and $\varepsilon$ in function $f$ exists and they are location, scale and skewness parameters, respectively. There are other parameters in the considered $f$, however we are not interested in other parameters.

In our proposal, our aim is to estimate the parameters $\theta$, $\sigma$ and $\varepsilon$ for the random sample $X_n=\{x_1,x_2,...,x_n\}$. The random sample is supposed to be asymmetrically distributed. Owing to the fact that the probability density function is not known, using the maximum likelihood estimation (MLE) method is not possible. In such a case, the function $Q$ given in the following form is proposed to estimate the parameters interested.

\begin{eqnarray}\label{minasymobjectivefunction}
 Q(\theta,\sigma,\varepsilon;X_n)   &=&  \sum_{i=1}^n \rho_{ESH} \bigg(\frac{x_i - \theta}{\sigma(1-sign(x_i-\theta)\varepsilon)} \bigg)+ n log(\sigma)\nonumber \\
   &+& \sum_{i=1}^n log(1-sign(x_i-\theta)\varepsilon)
\end{eqnarray}

The function in equation \eqref{minasymobjectivefunction} is minimized. To get the estimators of parameters $\theta$, $\sigma$ and $\varepsilon$, we will take the derivative of parameters interested.

Let $u_i$ be $\frac{x_i - \theta}{\sigma(1-sign(x_i - \theta)\varepsilon)}$. Then, the derivative of $Q(\theta,\sigma,\varepsilon;X_n)$ with respect to $\theta$ and setting it to zero will produce the following equation.
\begin{equation}\label{minasymobjectivefunction.theta}
  \frac{\partial}{\partial \theta}  Q(\theta,\sigma,\varepsilon;X_n) =  \sum_{i=1}^n \psi_\theta \bigg(\frac{x_i - \theta}{\sigma(1-sign(x_i-\theta)\varepsilon)}\bigg)\frac{-1}{\sigma(1-sign(x_i - \theta)\varepsilon)}=0.
\end{equation}
\noindent The weight function is defined to be $w(u_i)=\psi_\theta(u_i)/u_i$. Then, the M-estimator of location parameter $\theta$ will be
\begin{equation}\label{thsimulsteps}
  \hat{\theta}=\sum_{i=1}^n w_i\frac{x_i}{(\hat{\sigma}(1-sign(x_i-\hat{\theta})\hat{\varepsilon}))^2}/\sum_{i=1}^n w_i\frac{1}{(\hat{\sigma}(1-sign(x_i-\hat{\theta})\hat{\varepsilon}))^2}.
\end{equation}
\noindent where $w_i=w(\frac{x_i - \hat{\theta}}{\hat{\sigma}(1-sign(x_i - \hat{\theta})\hat{\varepsilon})})$. The derivative of $Q(\theta,\sigma,\varepsilon;X_n)$ with respect to $\sigma$ and setting it to zero will produce the following equation.
\begin{equation}\label{minasymobjectivefunction.sigma}
 \frac{\partial}{\partial \sigma}  Q(\theta,\sigma,\varepsilon;X_n) = -  \sum_{i=1}^n \psi_\sigma \bigg(\frac{x_i - \theta}{\sigma(1-sign(x_i-\theta)\varepsilon)}\bigg) \frac{x_i-\theta}{\sigma^2(1-sign(x_i - \theta)\varepsilon)}+\frac{n}{\sigma}=0.
\end{equation}
\noindent The weight function is defined to be $w(u_i)=\psi_\sigma(u_i)/u_i$. Then, the M-estimator of scale parameter $\sigma$ will be
\begin{equation}\label{sgsimulsteps}
  \hat{\sigma}^2=\frac{1}{n} \sum_{i=1}^n w_i \frac{(x_i - \hat{\theta})^2}{(1 - sign(x_i - \hat{\theta})\hat{\varepsilon})^2}.
\end{equation}
\noindent where $w_i=w(\frac{x_i -\hat{\theta}}{\hat{\sigma}(1-sign(x_i -\hat{\theta})\hat{\varepsilon})})$. The derivative of $Q(\theta,\sigma,\varepsilon;X_n)$ with respect to $\varepsilon$ and setting it to zero will produce the following equation.
\begin{eqnarray}\label{minasymobjectivefunction.epsilon}
\frac{\partial}{\partial \varepsilon} Q(\theta,\sigma,\varepsilon;X_n)  &=& \sum_{i=1}^n \psi_\varepsilon \bigg(\frac{x_i - \theta}{\sigma(1-sign(x_i-\theta)\varepsilon)} \bigg)\frac{(x_i -\theta) sign(x_i - \theta)}{\sigma(1-sign(x_i - \theta)\varepsilon)^2} \nonumber \\
   &-& \sum_{i=1}^{n}\frac{sign(x_i - \theta)}{(1-sign(x_i - \theta)\varepsilon)}=0.
\end{eqnarray}
\noindent  The weight function is defined to be $w(u_i)=\psi_\varepsilon(u_i)/u_i$. Then, the M-estimator of skewness parameter $\varepsilon$ will be
\begin{equation}\label{epsimulsteps}
  \hat{\varepsilon}=\sum_{i=1}^n\bigg[\frac{sign(x_i -\hat{\theta})}{(1-sign(x_i -\hat{\theta})\hat{\varepsilon})^2} -  w_i \frac{ (x_i - \hat{\theta})^2sign(x_i - \hat{\theta})}{\hat{\sigma}^2(1-sign(x_i-\hat{\theta})\hat{\varepsilon})^3}\bigg]/\sum_{i=1}^n\frac{1}{(1- sign(x_i - \hat{\theta})\hat{\varepsilon})^2}.
\end{equation}
\noindent The weight function is defined to be $w_i=w(\frac{x_i -\hat{\theta}}{\hat{\sigma}(1-sign(x_i -\hat{\theta})\hat{\varepsilon})})$.

The weight function of these parameters will be given as follow:
\begin{equation}\label{boundwithcepsilonWAM}
  w(u) = \left\{
                \begin{array}{ll}
                  \frac{c_1}{(1+\varepsilon)^2u}, & (-\infty,c_1]; \\
                  \frac{1}{(1+\varepsilon)^2}, & [c_1,0); \\
                  \frac{1}{(1-\varepsilon)^2}, & [0,c_2]; \\
                  \frac{c_2}{(1-\varepsilon)^2u}, & [c_2,\infty).
                \end{array}
              \right.
\end{equation}
\noindent As a result, the estimators of parameters $\theta$, $\sigma$  and $\varepsilon$  are gotten. The weight function in equation \eqref{boundwithcepsilonWAM} can give the different weights in data set that is negative and positive sides of axis. Thus, the estimators can model the asymmetry in the data set.

\subsection{The Computation Steps of Estimators}

The random sample is $X_n = \{x_1,x_2,...,x_n\}$ and $k \in \mathbb{N}^{+}$ is the iteration number. Then iterative reweighting algortihm (IRA) will be given in the following form:

\textbf{1. Step} $\theta^{(1)}$, $\sigma^{(1)}$ ve $\varepsilon^{(1)}$ are the initial values to start the algorithm.

\textbf{2. Step} The weight function $w$ in equation \eqref{boundwithcepsilonWAM} is computed by using the following form:
\begin{eqnarray*}
  u_i^{(k)} &=& \frac{x_i - \hat{\theta}^{(k)}}{\hat{\sigma}^{(k)}(1-sign(x_i - \hat{\theta}^{(k)})\hat{\varepsilon}^{(k)})}
\end{eqnarray*}

\textbf{3. Step} The estimated value of parameter $\theta$ is computed by
\begin{eqnarray*}
  \hat{\theta}^{(k+1)} &=& \sum_{i=1}^n w_i^{(k)} \frac{x_i}{(\hat{\sigma}^{{(k)}})^2(1-sign(x_i-\hat{\theta}^{(k)})\hat{\varepsilon}^{(k)})^2} \\
  && /\sum_{i=1}^n w_i^{(k)} \frac{1}{(\hat{\sigma}^{{(k)}})^2(1-sign(x_i-\hat{\theta}^{(k)})\hat{\varepsilon}^{(k)})^2}.
\end{eqnarray*}

\textbf{4. Step} The estimated value of parameter $\sigma$ is computed by
\begin{eqnarray*}
  (\hat{\sigma}^2)^{(k+1)} &=& \frac{1}{n}\sum_{i=1}^n w_i^{(k)}\frac{(x_i - \hat{\theta}^{(k+1)})^2}{(1 - sign(x_i - \hat{\theta}^{(k+1)})\hat{\varepsilon}^{(k)})^2}
\end{eqnarray*}

\textbf{5. Step} The estimated value of parameter $\varepsilon$ is computed by
\begin{eqnarray*}
  \hat{\varepsilon}^{(k+1)} &=&  \sum_{i=1}^n \bigg[\frac{sign(x_i -\hat{\theta}^{(k+1)})}{(1-sign(x_i -\hat{\theta})\hat{\varepsilon}^{(k)})^2} - w_i^{(k+1)}\frac{ (x_i - \hat{\theta}^{(k+1)})^2sign(x_i - \hat{\theta}^{(k+1)})}{(\hat{\sigma}^{(k+1)})^2(1-sign(x_i-\hat{\theta}^{(k+1)})\hat{\varepsilon}^{(k)})^3}\bigg] ~~\\
   && /\sum_{i=1}^n\frac{1}{(1- sign(x_i - \hat{\theta}^{(k+1)})\hat{\varepsilon}^{(k)})^2}
\end{eqnarray*}
\noindent where $w_i^{(k+1)}=w(u_i^{(k+1)})$. Then, the weight function $w$  in second step is computed by using the estimates $\hat{\theta}^{(k+1)}$, $\hat{\sigma}^{(k+1)}$ and $\hat{\varepsilon}^{(k)}$.

 \textbf{6. Step} If the norm of vector $(\hat{\theta}^{(k+1)} - \hat{\theta}^{(k)},\hat{\sigma}^{(k+1)}-\hat{\sigma}^{(k)}, \hat{\varepsilon}^{(k+1)}-\hat{\varepsilon}^{(k)})^T$ is bigger than the prescribed value $\epsilon>0$, the steps are repeated until  the prescribed value $\epsilon>0$ is guaranteed. Finally, the values at last steps are assigned to be estimates of parameters.

\section{Robustness Properties of Estimators}
In this section, the robustness properties of estimators of parameters $\theta,\sigma$ and $\varepsilon$  will be examined. In this context, the influence function that is an indicator for the local robustness and gross error sensitivity that is an indicator of global robustness are considered for the estimators of $\theta,\sigma$ and $\varepsilon$.

\subsection{The Influence Function of Estimators}

The function $\rho$ in equation \eqref{boundwithcepsilonrho} is used to get the following functions
$\psi_\theta(x) = \frac{\partial}{\partial \theta} \rho_{ESH}(\frac{x-\theta}{\sigma(1-sign(x)\varepsilon)}), \psi_\sigma(x)= \frac{\partial}{\partial \sigma} \rho_{ESH}(\frac{x-\theta}{\sigma(1-sign(x)\varepsilon)})$ and $\psi_\varepsilon(x)= \frac{\partial}{\partial \varepsilon} \rho_{ESH}(\frac{x-\theta}{\sigma(1-sign(x)\varepsilon)})$. For $\theta=0$ and $\sigma=1$, these functions are given by

\begin{equation}\label{psitheta01}
 \psi_\theta(x)=\left\{
                  \begin{array}{ll}
                   \frac{-c_1}{(1+\varepsilon)^3} , & (-\infty,c_1(1+\varepsilon)] ; \\
                   \frac{-x}{(1+\varepsilon)^4} , & [c_1(1+\varepsilon),0) ; \\
                    \frac{-x}{(1-\varepsilon)^4}, & [0,c_2(1-\varepsilon)]; \\
                     \frac{-c_2}{(1-\varepsilon)^3}, & [0+c_2(1-\varepsilon),\infty).
                  \end{array}
                \right.
\end{equation}
\begin{equation}\label{psisigma01}
\psi_\sigma(x)=\left\{
                  \begin{array}{ll}
                   \frac{-c_1x}{(1+\varepsilon)^3} , &(-\infty,c_1(1+\varepsilon)]  ; \\
                   \frac{-x^2}{(1+\varepsilon)^4} , & [c_1(1+\varepsilon),0); \\
                    \frac{-x^2}{(1-\varepsilon)^4}, & [0,c_2(1-\varepsilon)]; \\
                     \frac{-c_2x}{(1-\varepsilon)^3}, &[c_2(1-\varepsilon),\infty).
                  \end{array}
                \right.
\end{equation}
\begin{equation}\label{psiepsilon01}
 \psi_\varepsilon(x)=\left\{
                  \begin{array}{ll}
                   \frac{-3c_1x}{(1+\varepsilon)^4}+\frac{c_1^2}{(1+\varepsilon)^3} , &(-\infty,c_1(1+\varepsilon)]   ; \\
                   \frac{-2x^2}{(1+\varepsilon)^5} , & [c_1(1+\varepsilon),0) ; \\
                    \frac{2x^2}{(1-\varepsilon)^5}, & [0,c_2(1-\varepsilon)]; \\
                     \frac{3c_2x}{(1-\varepsilon)^4}-\frac{c_2^2}{(1-\varepsilon)^3}, &[c_2(1-\varepsilon),\infty).
                  \end{array}
                \right. 
\end{equation}

$\underset{x \rightarrow -\infty}{lim} \psi_\theta(x)= \frac{-c_1}{(1+\varepsilon)^3}<\infty$, however  $\underset{x \rightarrow -\infty}{lim}  \psi_\sigma(x) =  -\infty$ and $\underset{x \rightarrow -\infty}{lim}  \psi_\varepsilon(x) =-\infty$. Then, they are not finite. When the parameters $\sigma$ and $\varepsilon$ are known, the influence function of estimator of parameter $\theta$ is finite. However, the influence function of estimators of three parameters are not finite at the same time as it is proved by the tools given by the following forms:

The influence function of estimators of three parameters is
\begin{equation}\label{IFBExpectation}
IF(x;\hat{\theta},\hat{\sigma},\hat{\varepsilon})=- B^{-1}\Psi(x)
\end{equation}
\noindent where  $\Psi(x)=(\psi_\theta(x),\psi_\sigma(x),\psi_\varepsilon(x))^T$ and
\begin{equation}\label{Condition7MatrixB}
 B = \begin{bmatrix}
       E_{ESN}[\frac{\partial}{\partial\theta} \psi_\theta(X)] & E_{ESN}[\frac{\partial}{\partial \sigma}\psi_\theta(X)] & E_{ESN}[\frac{\partial}{\partial \varepsilon}\psi_\theta(X)]           \\[0.3em]
        &  E_{ESN}[\frac{\partial}{\partial\sigma}\psi_\sigma(X)]           & E_{ESN}[\frac{\partial}{\partial \varepsilon}\psi_\sigma(X)] \\[0.3em]
                  &  &  E_{ESN}[\frac{\partial}{\partial\varepsilon}\psi_\varepsilon(X)]
     \end{bmatrix}
\end{equation}

$ESN$ shows that the underlying distribution is taken as $ESN$ to get the integral values. It should be noted that $det(B) \neq 0$. Then, the matrix $B^{-1}$ exists. Thus, the influence function of estimators of three parameters exists.

The equation \eqref{IFBExpectation} can be rewritten as the following form:

\begin{equation}\label{comIFESEPeqns}
~~~~~~~ IF(x;\hat{\theta},\hat{\sigma},\hat{\varepsilon})=  \begin{bmatrix}
 T_{11}\psi_\theta(x) + T_{12}\psi_\sigma(x) +  T_{13} \psi_\varepsilon(x)  \\
 T_{21} \psi_\theta(x) +  T_{22}\psi_\sigma(x) +   T_{23} \psi_\varepsilon(x)  \\
  T_{31}\psi_\theta(x) +  T_{32} \psi_\sigma(x) +  T_{33}\psi_\varepsilon(x)
  \end{bmatrix} = \begin{bmatrix} IF_{1}(x;\hat{\theta},\hat{\sigma},\hat{\varepsilon}) \\ IF_{2}(x;\hat{\theta},\hat{\sigma},\hat{\varepsilon}) \\ IF_{3}(x;\hat{\theta},\hat{\sigma},\hat{\varepsilon})   \end{bmatrix}
\end{equation}

\noindent where $ T_{ij}$ represents the row $i.$ and column $j.$ of  matrix $B^{-1}$ ($i,j=1,2,3$). Here,  the components $IF_1,IF_2$ and $IF_3$ of the influence function $(IF)$ are not finite, because $\psi_\sigma(x)$ and $\psi_\varepsilon(x)$ are not finite. Thus, the influence function of the estimators is not finite. It is known that the norm of influence function is defined to be the gross error sensitivity. Then, the gross error sensitivity is given by the following form:
\begin{equation}\label{GESAsymmetricM}
GES(\hat{\theta},\hat{\sigma},\hat{\varepsilon},\rho_{ESH})=\{(IF_{1})^2+(IF_{2})^2+(IF_{3})^2\}^{1/2}.
\end{equation}
\noindent   The components $IF_1,IF_2$ and $IF_3$  are not finite. Then, $GES(\hat{\theta},\hat{\sigma},\hat{\varepsilon},\rho_{ESH})$ will not be finite.

\subsection{Breakdown Point of Estimator for Location Parameter}

\cite{Hub84} and \cite{ZhangLi98} give the assumption for the breakdown properties of location M-estimator. The convexity in asymmetric case is satisfied. Then, these assumptions given below can be used to satisfy the M-estimator generated from the asymmetric objective function.
\begin{enumerate}
\item $\rho(0)=0$  \cite{Hub84} and \cite{ZhangLi98}.
 \item $\underset{|u| \rightarrow \infty}{lim} \rho (u) = \infty$ \cite{Hub84} and \cite{ZhangLi98}.
\item $\psi(u) = \frac{d}{du}\rho(u)$  is continuous for every point of $u$.  \cite{Hub84}.
\item  Let $u_0$ exist when $\psi(u)$ is nondecreasing for $0<u \leq u_0$ and nonincreasing for $u_0 < u < \infty$ for monotone $\psi$ functions  \cite{ZhangLi98}.
\end{enumerate}
Let us check these assumptions for asymmetric objective and influence functions.
\begin{enumerate}
   \item $\rho_{ESH}(0)=0,$
   \item $\underset{|u| \rightarrow \infty}{lim} \rho_{ESH} (u) = \infty,$
   \item The function $\psi$ is continuous for every point of $u$.
   \item The function $\psi(u)$ is increasing for the point $c_2$ given arbitrarily for $[0, c_2]$  and it is constant for $(c_2, \infty)$
\end{enumerate}
These assumptions given above are satisfied. Then, the location estimator obtained from the asymmetric objective function  $\rho_{ESH}$ has a global robustness that is breakdown point. The value of breakdown point is $1/2$.

\section{Asymptotic Properties}
The asymptotic properties that are consistency and asymptotic normality of estimators for the parameters $\theta,\sigma$ ve $\varepsilon$ will be examined in this section. The function $Q$ in equation \eqref{minasymobjectivefunction} is used to show the asymptotic property. The following equations can be obtained after taking the derivatives with respect to parameters and setting them to zero.
The explicit forms of the following equations were given by the equations  \eqref{minasymobjectivefunction.theta},\eqref{minasymobjectivefunction.sigma} and \eqref{minasymobjectivefunction.epsilon}.

\begin{eqnarray}
  \sum_{i=1}^{n} \psi_\theta(x_i;\theta,\sigma,\varepsilon) &=& 0 \\
  \sum_{i=1}^{n}\psi_\sigma(x_i;\theta,\sigma,\varepsilon)  &=&0  \\
  \sum_{i=1}^{n} \psi_\varepsilon(x_i;\theta,\sigma,\varepsilon)  &=&0
\end{eqnarray}

The simultaneous estimations of the parameters $\theta$, $\sigma$ and $\varepsilon$
\begin{equation}\label{esanlipsi}
  \sum_{i=1}^{n}  \Psi(x_i;\hat{\theta},\hat{\sigma},\hat{\varepsilon})=\mathbf{0}.
\end{equation}
\noindent where $ \Psi=(\psi_\theta,\psi_\sigma,\psi_\varepsilon)$. The approach given by \cite{Scholz65}  is adapted into the asymmetric M-estimation. Then, there is an one solution of equation  \eqref{esanlipsi} at least.

Suppose that $\hat{\theta}$ exists for each of $\hat{\sigma}$. Then,

\begin{equation}\label{psitheta}
  \sum_{i=1}^{n}  \psi_\theta(x_i;\theta,\sigma,\varepsilon)=0.
\end{equation}
\noindent the location estimation $\hat{\theta}$
\begin{equation*}
  \underset{1 \leq i \leq n~~~~~~}{min(x_i)} \leq \hat{\theta} \leq \underset{1 \leq i \leq n~~~~~~}{max(x_i)}
\end{equation*}
\noindent is satisfied. Thus, at least one solution can exist for the location estimation. When the $\hat{\sigma}$ changes from $0$ to $\infty$. The term
 \begin{equation}\label{psisigma}
  \sum_{i=1}^{n}  \psi_\sigma(x_i;\theta,\sigma,\varepsilon)=0.
\end{equation}
\noindent changes from $sup\{\psi_\sigma(x_i;\theta,\sigma,\varepsilon): x \in \mathbb{R}\}$ to $0$. The estimation of parameter $\varepsilon$ will be solution of the following equation.
 \begin{equation}\label{psivarepsilon}
  \sum_{i=1}^{n} \psi_\varepsilon(x_i;\theta,\sigma,\varepsilon)=0.
\end{equation}
 The solution is in interval $(-1,1)$, because $\varepsilon \in (-1,1)$. Thus, at least one solution of equation \eqref{esanlipsi} can be given.

To guarantee the uniqueness of  solution, the following two conditions are satisfied \cite{FujiHer00}:
\begin{enumerate}
  \item The function $\rho_{ESH}$ is differentiable.
  \item The Jacobian of equation \eqref{esanlipsi} exists and upper-left corner principal minors of matrix is non-zero.
\end{enumerate}
We will examine whether the conditions are satisfied.

\begin{enumerate}
  \item The proposed function $\rho_{ESH}$ is differentiable at each points of the interval $[c_1,c_2]$. However, the derivative of $\rho_{ESH}$ is zero at the intervals $(-\infty,c_1)$ and $(c_2,\infty)$.
  \item   To construct the equation \eqref{esanlipsi}
   \begin{equation*}
 \lambda(\tau)=E_{ESN}  \Psi(X,\tau), ~~~~ \tau =(\theta,\sigma,\varepsilon)
\end{equation*}
\noindent then
\begin{equation*}
B_{jk}=\frac{\partial\lambda_j}{\partial\tau_k}, ~~~~~~~ j,k=1,2,3
\end{equation*}
exists each term of matrix $B$
\end{enumerate}
Constructing the matrix $B$ in second condition was proposed by \cite{Hub81,Maretal06}.

This matrix will be given when the asymptotic normality of estimators. Left-upper three corners of the matrix are given as follows:

$K_1=|B_{11}|$,~~
$K_2= \begin{vmatrix}
 B_{11}& B_{12}   \\
 B_{21}& B_{22} \\
 \end{vmatrix}$,~~
$K_3= \begin{vmatrix}
 B_{11}& B_{12}  & B_{13} \\
 B_{21}  & B_{22}& B_{23} \\
 B_{31}  & B_{32} &   B_{33}
 \end{vmatrix}$

\noindent $K_1,K_2$ and $K_3$ should be non-zero.   Then, the solutions  of equation \eqref{esanlipsi} has an unique solution, because the conditions $1.$ and $2.$ are satisfied.

Uniqueness of solution was guaranteed. Then, we will examine whether the estimators  $\hat{\theta},\hat{\sigma}$ ve $\hat{\varepsilon}$ are consistent. The convexity can hold when the asymmetric form of objective function is proposed. Then, we can use the assumptions considered by \cite{Haberman89}.
\begin{enumerate}
  \item $E_{ESN}[\rho_{ESH}(X)] < \infty$
  \item $E_{ESN}[\psi_{\theta}(X)] < \infty$, $E_F[\psi_{\sigma}(X)]< \infty$, $E_F[\psi_{\varepsilon}(X)] < \infty$
\end{enumerate}

The result of $E_{ESN}[\rho_{ESH}(X)]$ is given by
\begin{small}
\begin{eqnarray}\label{cond41th}
E[\rho(X)]&=&\frac{-c_1}{\sqrt{2\pi}}\cdot\Gamma(1,\frac{c_1^2}{2(1+\varepsilon)^2})-\frac{c_1^2}{4\sqrt{\pi}(1+\varepsilon)}\cdot\Gamma(\frac{1}{2},\frac{c_1^2}{2(1+\varepsilon)^2}) +\frac{(1+\varepsilon)}{2\sqrt{\pi}}\\ \nonumber
&&\cdot\gamma(\frac{3}{2},\frac{c_1^2}{2(1+\varepsilon)^2}) +\frac{(1-\varepsilon)}{2\sqrt{\pi}}\cdot\gamma(\frac{3}{2},\frac{c_2^2}{2(1-\varepsilon)^2})+\frac{c_2}{\sqrt{2\pi}}
 \\ \nonumber
&&\cdot\Gamma(1,\frac{c_2^2}{2(1-\varepsilon)^2})-\frac{c_2^2}{4\sqrt{\pi}(1-\varepsilon)}\cdot\Gamma(\frac{1}{2},\frac{c_2^2}{2(1-\varepsilon)^2}). \nonumber
\end{eqnarray}
\end{small}

The results of $E_{ESN}[\psi_{\theta}(X)]$, $E_{ESN}[\psi_{\sigma}(X)]$, $E_{ESN}[\psi_{\varepsilon}(X)]$ are given by

\begin{eqnarray}\label{cond42th}
E[\psi_\theta(X)]&=&\frac{-c_1}{\sigma(1+\varepsilon)^{2}2\sqrt{\pi}}\cdot\Gamma(\frac{1}{2},\frac{c_1^2}{2(1+\varepsilon)^2})
+\frac{1}{\sigma(1+\varepsilon)\sqrt{2\pi}}\cdot\gamma(1,\frac{c_1^2}{2(1+\varepsilon)^2})
 \\ \nonumber
&&-\frac{1}{\sigma(1-\varepsilon)\sqrt{2\pi}}\cdot\gamma(1,\frac{c_2^2}{2(1-\varepsilon)^2})
-\frac{c_2}{\sigma(1-\varepsilon)^{2}2\sqrt{\pi}}\cdot\Gamma(\frac{1}{2},\frac{c_2^2}{2(1-\varepsilon)^2}), \nonumber
\end{eqnarray}
\begin{eqnarray}\label{cond42sg}
E[\psi_\sigma(X)]&=&\frac{c_1}{\sigma\sqrt{2\pi}}\cdot\Gamma(1,\frac{c_1^2}{2(1+\varepsilon)^2})-\frac{(1+\varepsilon)}{\sigma\sqrt{\pi}}\cdot\gamma(\frac{3}{2},\frac{c_1^2}{2(1+\varepsilon)^2})
\\ \nonumber
&&-\frac{(1-\varepsilon)}{\sigma\sqrt{\pi}}\cdot\gamma(\frac{3}{2},\frac{c_2^2}{2(1-\varepsilon)^2})-\frac{c_2}{\sigma\sqrt{2\pi}}\cdot\Gamma(1,\frac{c_2^2}{2(1-\varepsilon)^2}), \end{eqnarray}
\begin{eqnarray}\label{cond42ep}
E[\psi_\varepsilon(X)]&=&\frac{3c_1}{(1+\varepsilon)\sqrt{2\pi}}\cdot\Gamma(1,\frac{c_1^2}{2(1+\varepsilon)^2})+\frac{c_1^2}{(1+\varepsilon)^{2}2\sqrt{\pi}}\Gamma(\frac{1}{2},\frac{c_1^2}{2(1+\varepsilon)^2})
 \\ \nonumber
&&-\frac{2}{\sqrt{\pi}}\cdot\gamma(\frac{3}{2},\frac{c_1^2}{2(1+\varepsilon)^2})+\frac{2}{\sqrt{\pi}}\cdot\gamma(\frac{3}{2},\frac{c_2^2}{2(1-\varepsilon)^2})
 \\ \nonumber
&&+\frac{3c_2}{(1-\varepsilon)\sqrt{2\pi}}\Gamma(1,\frac{c_2^2}{2(1-\varepsilon)^2})-\frac{c_2^2}{(1-\varepsilon)^{2}2\sqrt{\pi}}\Gamma(\frac{1}{2},\frac{c_2^2}{2(1-\varepsilon)^2}). \end{eqnarray}

These results will be finite when the $c_1$, $c_2$, $\sigma$ and $\varepsilon$ are finite. Then, the conditions $1.$ and $2.$ were satisfied. The consistency of estimators obtained simultaneously was been examined. Then, it can examine the asymptotic normality of estimators. Since the estimators are not explicit form, the Taylor expansion of influence function will be consider as follow:

\begin{equation}\label{aaaVector4}
 \Psi(x_i,\hat{\tau})=\Psi(x_i,\tau)+(\hat{\tau}-\tau)\dot{\Psi}(x_i,\tau)+\mathbf{R}_n^{*}
\end{equation}
the summation of both sides are taken and multiplied by $1/n$. Then,
\begin{equation*}
 \mathbf{0}=\frac{1}{n}  \sum_{i=1}^n  \Psi(x_i,\tau)+(\hat{\tau}-\tau)\frac{1}{n}  \sum_{i=1}^n  \dot{\Psi}(x_i,\tau)+\mathbf{R}_n
\end{equation*}
\noindent was gotten, where $\dot{\Psi}(x_i,\tau) = \frac{\partial \Psi(x_i,\tau)}{\partial \tau^T}$. If $[\frac{1}{n}  \sum_{i=1}^n \dot{\Psi}(x_i,\tau)]^{-1}$  exists,
\begin{eqnarray}
 -(\hat{\tau} - \tau)\frac{1}{n} \sum_{i=1}^n \dot{\Psi}(x_i,\tau)  &=& \frac{1}{n} \sum_{i=1}^n \Psi(x_i,\tau)+\mathbf{R}_n \\ \nonumber
   \sqrt{n}(\hat{\tau} - \tau)  &=& B_n^{-1} \sqrt{n} A_n +\sqrt{n}\mathbf{R}_n   \nonumber
\end{eqnarray}
\noindent where $\sqrt{n}\mathbf{R}_n  \overset{P}{\longrightarrow}  \mathbf{0}$. Under the regularity of conditions, when $n \rightarrow \infty$, the weak of large numbers
\begin{equation}\label{convergenceinp4}
 B_n = \frac{1}{n} \sum_{i=1}^n (-\dot{\Psi}(x_i,\tau)) \overset{P}{\longrightarrow} E[-\dot{\Psi}(X,\tau)] = B
\end{equation}
\noindent can be obtained. By means of central limit theorem,
\begin{equation}\label{mlt}
  \sqrt{n}A_n  \overset{D}{\longrightarrow} N_3(0,A ), ~~~~A =E[\Psi(X,\tau)\Psi(X,\tau)^{T}]
\end{equation}
\noindent can be obtained. Here, $\Psi$ is three-dimensional. Thus, by means of the Slutsky's multivariate lemma,
\begin{equation}\label{AMNSlut4}
 \sqrt{n}(\hat{\tau}-\tau) \overset{D}{\longrightarrow} N_3(0,B^{-1}A(B^{T})^{-1})
\end{equation}
\noindent $N_3$ shows the $3-$dimensional asymptotic normally distributed.

Then,  the matrices $A$ and $B$ exist and the inverse of matrix $B$ exists as well. These matrices are obtained when the underlying distribution is chosen as $ESN$.

$A= \begin{bmatrix}
E[\psi_\theta^2(X)]& E[\psi_\sigma(X)\psi_\theta(X)]  & E[\psi_\varepsilon(X)\psi_\theta(X)] \\
   &E[\psi_\sigma^2(X)] & E[\psi_\varepsilon(X)\psi_\sigma(X)] \\
   &  &   E[\psi_\varepsilon^2(X)]
 \end{bmatrix}$
 
The elements of matrix $A$ are

\begin{eqnarray*}
E[\psi_\theta^2(X)]&=&\frac{c_1^2}{\sigma^2(1+\varepsilon)^{5}2\sqrt{\pi}}\cdot\Gamma(\frac{1}{2},\frac{c_1^2}{2(1+\varepsilon)^2})
+\frac{1}{\sigma^2(1+\varepsilon)^{3}\sqrt{\pi}}\cdot\gamma(\frac{3}{2},\frac{c_{1}^2}{2(1+\varepsilon)^2})
\nonumber \\
&&+\frac{1}{\sigma^2(1-\varepsilon)^{3}\sqrt{\pi}}\cdot\gamma(\frac{3}{2},\frac{c_{2}^2}{2(1-\varepsilon)^2})
+\frac{c_2^2}{\sigma^2(1-\varepsilon)^{5}2\sqrt{\pi}}\cdot\Gamma(\frac{1}{2},\frac{c_2^2}{2(1-\varepsilon)^2}),\\
E[\psi_\sigma(X) \psi_\theta(X)]&=&\frac{-c_1^2}{\sigma^2(1+\varepsilon)^3\sqrt{2\pi}}\cdot\Gamma(1,\frac{c_1^2}{2(1+\varepsilon)^2})
 -\frac{2}{\sigma^2(1+\varepsilon)\sqrt{2\pi}}\cdot\gamma(2,\frac{c_1^2}{2(1+\varepsilon)^2}) \\
 &&
 +\frac{2}{\sigma^2(1-\varepsilon)\sqrt{2\pi}}\cdot\gamma(2,\frac{c_2^2}{2(1-\varepsilon)^2})
 +\frac{c_2^2}{\sigma^2(1-\varepsilon)^3\sqrt{2\pi}}\cdot\Gamma(1,\frac{c_2^2}{2(1-\varepsilon)^2}), \\
\end{eqnarray*} 
\begin{eqnarray*}
E[\psi_\varepsilon(X) \psi_\theta(X)]&=&\frac{-3c_1^2}{\sigma(1+\varepsilon)^4\sqrt{2\pi}}\cdot\Gamma(1,\frac{c_1^2}{2(1+\varepsilon)^2})
-\frac{c_1^3}{\sigma(1+\varepsilon)^{5}2\sqrt{\pi}}\cdot\Gamma(\frac{1}{2},\frac{c_1^2}{2(1+\varepsilon)^2})
\nonumber \\
&&-\frac{2\sqrt{2}}{\sigma(1+\varepsilon)^2\sqrt{\pi}}\cdot\gamma(2,\frac{c_1^2}{2(1+\varepsilon)^2})
-\frac{2\sqrt{2}}{\sigma(1-\varepsilon)^2\sqrt{\pi}}\cdot\gamma(2,\frac{c_2^2}{2(1-\varepsilon)^2})
\nonumber \\
&&-\frac{3c_2^2}{\sigma(1-\varepsilon)^4\sqrt{2\pi}}\cdot\Gamma(1,\frac{c_2^2}{2(1-\varepsilon)^2})
 +\frac{c_2^3}{\sigma(1-\varepsilon)^{5}2\sqrt{\pi}}\cdot\Gamma(\frac{1}{2},\frac{c_2^2}{2(1-\varepsilon)^2}),\\
E[\psi_\sigma^2(X)]&=&\frac{c_1^2}{\sigma^2(1+\varepsilon)\sqrt{\pi}}\cdot\Gamma(\frac{3}{2},\frac{c_1^2}{2(1+\varepsilon)^2})+\frac{2(1+\varepsilon)}{\sigma^{2}\sqrt{\pi}}\cdot\gamma(\frac{5}{2},\frac{c_1^2}{2(1+\varepsilon)^2})
\nonumber \\
&&+\frac{2(1-\varepsilon)}{\sigma^{2}\sqrt{\pi}}\cdot\gamma(\frac{5}{2},\frac{c_2^2}{2(1-\varepsilon)^2})+\frac{c_2^2}{\sigma^2(1-\varepsilon)\sqrt{\pi}}\cdot\Gamma(\frac{3}{2},\frac{c_2^2}{2(1-\varepsilon)^2}), \\
E[\psi_\varepsilon(X)\psi_\sigma(X)]&=&\frac{3c_1^2}{\sigma(1+\varepsilon)^{2}\sqrt{\pi}}\cdot\Gamma(\frac{3}{2},\frac{c_1^2}{2(1+\varepsilon)^2})
+\frac{c_1^3}{\sigma(1+\varepsilon)^3\sqrt{2\pi}}\cdot\Gamma(1,\frac{c_1^2}{2(1+\varepsilon)^2})
\nonumber \\
&&+\frac{4}{\sigma\sqrt{\pi}}\cdot\gamma(\frac{5}{2},\frac{c_1^2}{2(1+\varepsilon)^2})-\frac{4}{\sigma\sqrt{\pi}}\cdot\gamma(\frac{5}{2},\frac{c_2^2}{2(1-\varepsilon)^2})
\nonumber \\
&&-\frac{3c_2^2}{\sigma(1-\varepsilon)^{2}\sqrt{\pi}}\Gamma(\frac{3}{2},\frac{c_2^2}{2(1-\varepsilon)^2})
+\frac{c_2^3}{\sigma(1-\varepsilon)^3\sqrt{2\pi}}\Gamma(1,\frac{c_2^2}{2(1-\varepsilon)^2}),\\
E[\psi_\varepsilon^2(X)]&=&\frac{9c_{1}^2}{(1+\varepsilon)^{3}\sqrt{\pi}}\cdot\Gamma(\frac{3}{2},\frac{c_{1}^2}{2(1+\varepsilon)^2})+\frac{3\sqrt{2}c_{1}^3}{(1+\varepsilon)^{4}\sqrt{\pi}}\cdot\Gamma(1,\frac{c_{1}^2}{2(1+\varepsilon)^2})
\nonumber \\
&&+\frac{c_{1}^4}{(1+\varepsilon)^{5}2\sqrt{\pi}}\cdot\Gamma(\frac{1}{2},\frac{c_{1}^2}{2(1+\varepsilon)^2})+\frac{8}{(1+\varepsilon)\sqrt{\pi}}\cdot\gamma(\frac{5}{2},\frac{c_{1}^2}{2(1+\varepsilon)^2})
\nonumber \\
&&+\frac{8}{(1-\varepsilon)\sqrt{\pi}}\cdot\gamma(\frac{5}{2},\frac{c_{2}^2}{2(1-\varepsilon)^2})+\frac{9c_{2}^2}{(1-\varepsilon)^{3}\sqrt{\pi}}\cdot\Gamma(\frac{3}{2},\frac{c_{2}^2}{2(1-\varepsilon)^2})
\nonumber \\
&&-\frac{3\sqrt{2}c_{2}^3}{(1-\varepsilon)^{4}\sqrt{\pi}}\cdot\Gamma(1,\frac{c_{2}^2}{2(1-\varepsilon)^2})+\frac{c_{2}^4}{(1-\varepsilon)^{5}2\sqrt{\pi}}\cdot\Gamma(\frac{1}{2},\frac{c_{2}^2}{2(1-\varepsilon)^2}).
\end{eqnarray*}
 
 $B=- \begin{bmatrix}
E[\frac{\partial}{\partial \theta}\psi_\theta(X) ]&  E[\frac{\partial}{\partial \sigma }\psi_\theta(X) ] & E[\frac{\partial}{\partial \varepsilon }\psi_\theta(X) ]\\
   &E[\frac{\partial}{\partial \sigma}\psi_\sigma(X) ]& E[\frac{\partial}{\partial \varepsilon}\psi_\sigma(X) ]\\
   &  &  E[\frac{\partial}{\partial \varepsilon}\psi_\varepsilon(X)]
 \end{bmatrix}$
 
 The elements of matrix $B$ are
 
 \begin{eqnarray*}
~~~~~~~~~~~~~~ E[\frac{\partial}{\partial\theta} \psi_\theta(X)]&=&\frac{1}{(1+\varepsilon)^{3}\sigma^{2}2\sqrt{\pi}}\cdot\gamma(\frac{1}{2},\frac{c_1^2}{2(1+\varepsilon)^2})
+\frac{1}{(1-\varepsilon)^{3}\sigma^{2}2\sqrt{\pi}}\cdot\gamma(\frac{1}{2},\frac{c_2^2}{2(1-\varepsilon)^2}),\\
E[\frac{\partial}{\partial \sigma}\psi_\theta(X)]&=&\frac{c_1}{\sigma^{2}(1+\varepsilon)^{2}2\sqrt{\pi}}\cdot\Gamma(\frac{1}{2},\frac{c_1^2}{2(1+\varepsilon)^2})
-\frac{\sqrt{2}}{\sigma^{2}(1+\varepsilon)\sqrt{\pi}}\cdot\gamma(1,\frac{c_1^2}{2(1+\varepsilon)^2})
\nonumber \\
&&+\frac{\sqrt{2}}{\sigma^{2}(1-\varepsilon)\sqrt{\pi}}\gamma(1,\frac{c_2^2}{2(1-\varepsilon)^2})
+\frac{c_2}{\sigma^{2}(1-\varepsilon)^{2}2\sqrt{\pi}}\cdot\Gamma(\frac{1}{2},\frac{c_2^2}{2(1-\varepsilon)^2}),\\
E[\frac{\partial}{\partial \varepsilon}\psi_\theta(X)]&=&\frac{3c_1}{\sigma(1+\varepsilon)^{3}2\sqrt{\pi}}\cdot\Gamma(\frac{1}{2},\frac{c_1^2}{2(1+\varepsilon)^2})
-\frac{2\sqrt{2}}{\sigma(1+\varepsilon)^2\sqrt{\pi}}\gamma(1,\frac{c_1^2}{2(1+\varepsilon)^2})
\nonumber \\
&&-\frac{2\sqrt{2}}{\sigma(1-\varepsilon)^2\sqrt{\pi}}\gamma(1,\frac{c_2^2}{2(1-\varepsilon)^2})
-\frac{3c_2}{\sigma(1-\varepsilon)^{3}2\sqrt{\pi}}\cdot\Gamma(\frac{1}{2},\frac{c_2^2}{2(1-\varepsilon)^2}),\\
E[\frac{\partial}{\partial\sigma}\psi_\sigma(X)]&=&-\frac{\sqrt{2}c_1}{\sigma^2\sqrt{\pi}}\cdot\Gamma(1,\frac{c_1^2}{2(1+\varepsilon)^2})+\frac{3(1+\varepsilon)}{\sigma^{2}\sqrt{\pi}}\cdot\gamma(\frac{3}{2},\frac{c_1^2}{2(1+\varepsilon)^2})
\nonumber \\
&&+\frac{3(1-\varepsilon)}{\sigma^{2}\sqrt{\pi}}\cdot\gamma(\frac{3}{2},\frac{c_2^2}{2(1-\varepsilon)^2})+\frac{\sqrt{2}c_2}{\sigma^2\sqrt{\pi}}\cdot\Gamma(1,\frac{c_2^2}{2(1-\varepsilon)^2}),\\
E[\frac{\partial}{\partial \varepsilon}\psi_\sigma(X)]&=&\frac{-3c_1}{\sigma(1+\varepsilon)\sqrt{2\pi}}\cdot\Gamma(1,\frac{c_1^2}{2(1+\varepsilon)^2})
+\frac{4}{\sigma\sqrt{\pi}}\cdot\gamma(\frac{3}{2},\frac{c_1^2}{2(1+\varepsilon)^2})
\nonumber \\
&&
-\frac{4}{\sigma\sqrt{\pi}}\cdot\gamma(\frac{3}{2},\frac{c_2^2}{2(1-\varepsilon)^2})
-\frac{3c_2}{\sigma(1-\varepsilon)\sqrt{2\pi}}\cdot\Gamma(1,\frac{c_2^2}{2(1-\varepsilon)^2}),\\
E[\frac{\partial}{\partial\varepsilon}\psi_\varepsilon(X)]&=&\frac{-6\sqrt{2}c_1}{(1+\varepsilon)^2\sqrt{\pi}}\cdot\Gamma(1,\frac{c_{1}^2}{2(1+\varepsilon)^2})-\frac{3c_{1}^2}{(1+\varepsilon)^{3}2\sqrt{\pi}}\cdot\Gamma(\frac{1}{2},\frac{c_{1}^2}{2(1+\varepsilon)^2})
\nonumber \\
&&+\frac{10}{(1+\varepsilon)\sqrt{\pi}}\cdot\gamma(\frac{3}{2},\frac{c_{1}^2}{2(1+\varepsilon)^2})+\frac{10}{(1-\varepsilon)\sqrt{\pi}}\cdot\gamma(\frac{3}{2},\frac{c_{2}^2}{2(1-\varepsilon)^2})
\nonumber \\
&&+\frac{6\sqrt{2}c_2}{(1-\varepsilon)^2\sqrt{\pi}}\Gamma(1,\frac{c_{2}^2}{2(1-\varepsilon)^2})-\frac{3c_{2}^2}{(1-\varepsilon)^{3}2\sqrt{\pi}}\Gamma(\frac{1}{2},\frac{c_{2}^2}{2(1-\varepsilon)^2}).
\end{eqnarray*}

 \begin{table}[!]
\caption{Asymptotic variance of estimators}
\centering
\begin{tabular}{cccccc}
\hline
                                                         & $Var(\hat{\tau})/n$                   & $n=30$ & $n=50$ & $n=100$ & $n=150$ \\ \hline \hline
                             & $Var(\hat{\theta})/n$                 & $0.190253 $    & $0.114152 $    & $0.057076$     & $0.038051 $     \\
$\varepsilon=-0.2$ & $Var(\hat{\sigma})/n$                & $0.018747 $    & $0.011248 $    & $0.005624 $     & $0.003749 $     \\
$c_1=-1.1,c_2=3.7$ & $Var(\hat{\varepsilon})/n$        & $0.021061$    & $0.012637 $    & $0.006318 $     & $0.004212 $     \\ \hline

                              & $Var(\hat{\theta})/n$                   & $0.059406 $    & $0.035644 $    & $0.017822$     &  $0.011881 $     \\
$\varepsilon=-0.5$ & $Var(\hat{\sigma})/n$                 & $0.022944 $    & $0.013767 $    & $0.006883 $     & $0.004589 $     \\
$c_1=-0.7,c_2=5.0$ & $Var(\hat{\varepsilon})/n$         & $0.016147 $    & $0.009688 $    & $0.004844 $     & $0.003229$     \\ \hline

                              & $Var(\hat{\theta})/n$                 & $0.010032 $    & $0.006019 $    & $0.003010 $     & $0.002006 $     \\
$\varepsilon=-0.8$ & $Var(\hat{\sigma})/n$                 & $0.035191 $    & $0.021114 $    & $0.010557 $     & $0.007038 $     \\
$c_1=-0.1,c_2=6.4$ & $Var(\hat{\varepsilon})/n$         & $0.023486 $    & $0.014091 $    & $0.007046 $     & $0.004697 $     \\ \hline
\end{tabular}
\label{ESHRaoCramerlowerbound}
\end{table}

\subsection{Simulation Study for Estimators of Location, Scale and Skewness Parameters}

To test the performance of asymmetric M-estimators, the contaminated version of ESN distribution will be considered as follow:
 $$0.90ESN(\theta=0,\sigma=1,\varepsilon=\varepsilon_0)+0.10ESL(\theta=0,\sigma=1,\varepsilon=\varepsilon_0)$$
$1000$ runs are performed. The sample sizes of each run are $30,50,100$ and $150$. The relative efficiencies of estimators are also computed.

\begin{eqnarray}\label{REMSERE}
  RE_{ESH}(\hat{\tau}) &=& \bigg(\frac{{MSE}_{ESH}(\hat{\tau})}{{MSE}_{ESN}(\hat{\tau})}\bigg)100 \\
  RE_{ESH}(\hat{\tau}) &=& \bigg(\frac{{MSE}_{ESH}(\hat{\tau})}{{MSE}_{ESL}(\hat{\tau})}\bigg)100  \\
  RE_{ESH}(\hat{\tau}) &=& \bigg(\frac{{MSE}_{ESH}(\hat{\tau})}{{MSE}_{ESt}(\hat{\tau})} \bigg)100 \\
  RE_{ESH}(\hat{\tau}) &=& \bigg(\frac{{MSE}_{ESH}(\hat{\tau})}{{MSE}_{H}(\hat{\tau})} \bigg)100
\end{eqnarray}

MSE is mean squared error obtained from simulation. ESH is epsilon-skew Huber M-estimator, ESN, ESL and ESt epsilon-skew normal, Laplace and $t$ distributions, respectively. H is symmetric Huber M-estimator. Three degrees of asymmetry are considered to test the asymmetric M-estimators. We make a comparison betwen the symmetric Huber M-estimator and asymmetric Huber M-estimator. In tables, maximum likelihood estimators of location, scale and skewness parameters of ESN, ESL and epsilon-skew t (ESt) considered by \cite{ArslanGenc09} distributions are given. The comparison of them with asymmetric M-estimator is also considered. Tables 2-\ref{ESHcomsen0835n} shows that asymmetric M-estimator (ESH) outperforms generally than the maximum likelihood and M-estimators when the data set has outliers.  The inital points of $\theta$, $\sigma$  and $\varepsilon$ to start the algorithm are median(x),  $MAD(x)=median(|x_i - median(x)|)$ and $0$, respectively. Here, $x=\{x_1,x_2,...,x_n\}$. Three degress of asymmetry are considered.

\begin{table}[t]
\centering
\label{ESHcomsen0235n}
\caption{Asymmetric M (ESH), ML and M Estimators ($\varepsilon=-0.2$):$c_1=   -1.10, c_2=    3.70$, $k=1.4$}
\scalebox{0.9}{
\begin{tabular}{cccccccccc}\hline
 &   & & $n=      30$ & & &  &$n=      50$  & &  \\ \hline\hline
& $\tau$ & $\hat{\tau}$ & ${Var}(\hat{\tau})$ & ${MSE}(\hat{\tau})$ & $RE$ & $\hat{\tau}$ & ${Var}(\hat{\tau})$ & ${MSE}(\hat{\tau})$ & $RE$ \\ \hline
\multicolumn{10}{c}{$ESH$} \\ \hline
$\theta$&$0.0$&$ -0.0357$&$  0.1645$&$  0.1658$&$100$&$ -0.0340$&$  0.0944$&$  0.0956$&$100$\\
$\sigma$&$1.0$&$  1.1180$&$  0.0633$&$  0.0773$&$100$&$  1.1566$&$  0.0348$&$  0.0594$&$100$\\
$\varepsilon$&$-0.2$&$ -0.1603$&$  0.0200$&$  0.0216$&$100$&$ -0.1804$&$  0.0153$&$  0.0160$&$100$\\
\hline\hline
\multicolumn{10}{c}{$ESN$} \\ \hline
$\theta$&$0.0$&$  0.3057$&$  0.0911$&$  0.1846$&$90$&$  0.3209$&$  0.0181$&$  0.1211$&$79$\\
$\sigma$&$1.0$&$  1.2466$&$  0.1105$&$  0.1714$&$45$&$  1.2999$&$  0.0696$&$  0.1595$&$37$\\
$\varepsilon$&$-0.2$&$ -0.1462$&$  0.0569$&$  0.0598$&$36$&$ -0.2200$&$  0.0454$&$  0.0455$&$35$\\
\hline\hline
\multicolumn{10}{c}{$ESL$} \\ \hline
$\theta$&$0.0$&$  0.1736$&$  0.1014$&$  0.1315$&$126$&$  0.1809$&$  0.0565$&$  0.0892$&$107$\\
$\sigma$&$1.0$&$  0.6316$&$  0.0137$&$  0.1494$&$52$&$  0.6499$&$  0.0074$&$  0.1300$&$46$\\
$\varepsilon$&$-0.2$&$ -0.1306$&$  0.0308$&$  0.0356$&$61$&$ -0.1263$&$  0.0172$&$  0.0226$&$71$\\
\hline\hline
\multicolumn{10}{c}{$ESt$} \\ \hline
$\theta$&$0.0$&$  0.0096$&$  0.2177$&$  0.2178$&$76$&$  0.0257$&$  0.1335$&$  0.1341$&$71$\\
$\sigma$&$1.0$&$  0.6850$&$  0.0145$&$  0.1137$&$68$&$  0.7136$&$  0.0081$&$  0.0902$&$66$\\
$\varepsilon$&$-0.2$&$ -0.2156$&$  0.0873$&$  0.0876$&$25$&$ -0.2007$&$  0.0510$&$  0.0510$&$31$\\
\hline
\multicolumn{10}{c}{$Huber$ $M$} \\ \hline
$\theta$&$0.0$&$  0.3342$&$  0.0714$&$  0.1831$&$91$&$  0.3640$&$  0.0388$&$  0.1714$&$56$\\
$\sigma$&$1.0$&$  1.1378$&$  0.0425$&$  0.0615$&$126$&$  1.0848$&$  0.0222$&$  0.0294$&$202$\\
\hline
 &   & & $n=     100$ & & &  &$n=     150$  & &  \\ \hline\hline
& $\tau$ & $\hat{\tau}$ & ${Var}(\hat{\tau})$ & ${MSE}(\hat{\tau})$ & $RE$ & $\hat{\tau}$ & ${Var}(\hat{\tau})$ & ${MSE}(\hat{\tau})$ & $RE$ \\ \hline
\multicolumn{10}{c}{$ESH$} \\ \hline
$\theta$&$0.0$&$ -0.0080$&$  0.0580$&$  0.0581$&$100$&$  0.0176$&$  0.0364$&$  0.0367$&$100$\\
$\sigma$&$1.0$&$  1.1555$&$  0.0197$&$  0.0439$&$100$&$  1.1702$&$  0.0106$&$  0.0396$&$100$\\
$\varepsilon$&$-0.2$&$ -0.1708$&$  0.0049$&$  0.0057$&$100$&$ -0.1891$&$  0.0011$&$  0.0012$&$100$\\
\hline\hline
\multicolumn{10}{c}{$ESN$} \\ \hline
$\theta$&$0.0$&$  0.1596$&$  0.0616$&$  0.0871$&$67$&$  0.1904$&$  0.0213$&$  0.0575$&$64$\\
$\sigma$&$1.0$&$  1.2967$&$  0.0421$&$  0.1331$&$33$&$  1.3220$&$  0.0274$&$  0.1311$&$30$\\
$\varepsilon$&$-0.2$&$ -0.1810$&$  0.0230$&$  0.0234$&$24$&$ -0.2346$&$  0.0094$&$  0.0106$&$11$\\
\hline\hline
\multicolumn{10}{c}{$ESL$} \\ \hline
$\theta$&$0.0$&$  0.1711$&$  0.0395$&$  0.0688$&$84$&$  0.1566$&$  0.0223$&$  0.0468$&$78$\\
$\sigma$&$1.0$&$  0.6492$&$  0.0043$&$  0.1274$&$34$&$  0.6531$&$  0.0024$&$  0.1228$&$32$\\
$\varepsilon$&$-0.2$&$ -0.1352$&$  0.0093$&$  0.0135$&$42$&$ -0.1321$&$  0.0066$&$  0.0112$&$11$\\
\hline\hline
\multicolumn{10}{c}{$ESt$} \\ \hline
$\theta$&$0.0$&$ -0.0206$&$  0.0741$&$  0.0746$&$78$&$  0.0004$&$  0.0446$&$  0.0446$&$82$\\
$\sigma$&$1.0$&$  0.7239$&$  0.0041$&$  0.0804$&$55$&$  0.7160$&$  0.0030$&$  0.0800$&$49$\\
$\varepsilon$&$-0.2$&$ -0.2270$&$  0.0252$&$  0.0259$&$22$&$ -0.2103$&$  0.0156$&$  0.0157$&$8$\\
\hline
\multicolumn{10}{c}{$Huber$ $M$} \\ \hline
$\theta$&$0.0$&$  0.3298$&$  0.0168$&$  0.1256$&$46$&$  0.3208$&$  0.0112$&$  0.1141$&$32$\\
$\sigma$&$1.0$&$  1.0469$&$  0.0184$&$  0.0206$&$213$&$  1.0640$&$  0.0095$&$  0.0136$&$291$\\
\hline
\end{tabular}
}
\end{table}

\begin{table}[t]
\centering
\label{ESHcomsen0535n}
\caption{Asymmetric M (ESH), ML and M Estimators ($\varepsilon=-0.5$):$c_1=   -0.70, c_2=    5.00$, $k=1.4$}
\scalebox{0.9}{
\begin{tabular}{cccccccccc}\hline
 &   & & $n=      30$ & & &  &$n=      50$  & &  \\ \hline\hline
& $\tau$ & $\hat{\tau}$ & ${Var}(\hat{\tau})$ & ${MSE}(\hat{\tau})$ & $RE$ & $\hat{\tau}$ & ${Var}(\hat{\tau})$ & ${MSE}(\hat{\tau})$ & $RE$ \\ \hline
\multicolumn{10}{c}{$ESH$} \\ \hline
$\theta$&$0.0$&$ -0.0828$&$  0.1534$&$  0.1602$&$100$&$ -0.0741$&$  0.1083$&$  0.1138$&$100$\\
$\sigma$&$1.0$&$  1.0458$&$  0.0546$&$  0.0567$&$100$&$  1.0274$&$  0.0369$&$  0.0377$&$100$\\
$\varepsilon$&$-0.5$&$ -0.4123$&$  0.0294$&$  0.0371$&$100$&$ -0.4415$&$  0.0137$&$  0.0171$&$100$\\
\hline\hline
\multicolumn{10}{c}{$ESN$} \\ \hline
$\theta$&$0.0$&$  0.4399$&$  0.1043$&$  0.2979$&$54$&$  0.3843$&$  0.0774$&$  0.2251$&$50$\\
$\sigma$&$1.0$&$  1.4055$&$  0.2091$&$  0.3735$&$15$&$  1.4629$&$  0.1720$&$  0.3863$&$10$\\
$\varepsilon$&$-0.5$&$ -0.3794$&$  0.0520$&$  0.0665$&$56$&$ -0.4288$&$  0.0341$&$  0.0392$&$44$\\
\hline\hline
\multicolumn{10}{c}{$ESL$} \\ \hline
$\theta$&$0.0$&$  0.4512$&$  0.1137$&$  0.3172$&$51$&$  0.3933$&$  0.0694$&$  0.2241$&$51$\\
$\sigma$&$1.0$&$  0.6787$&$  0.0200$&$  0.1233$&$46$&$  0.6933$&$  0.0118$&$  0.1059$&$36$\\
$\varepsilon$&$-0.5$&$ -0.2724$&$  0.0353$&$  0.0871$&$43$&$ -0.3064$&$  0.0207$&$  0.0582$&$29$\\
\hline\hline
\multicolumn{10}{c}{$ESt$} \\ \hline
$\theta$&$0.0$&$  0.1200$&$  0.2063$&$  0.2207$&$73$&$  0.0367$&$  0.1215$&$  0.1228$&$93$\\
$\sigma$&$1.0$&$  0.6934$&$  0.0145$&$  0.1085$&$52$&$  0.7172$&$  0.0068$&$  0.0868$&$43$\\
$\varepsilon$&$-0.5$&$ -0.4559$&$  0.0777$&$  0.0797$&$47$&$ -0.4947$&$  0.0443$&$  0.0443$&$39$\\
\hline
\multicolumn{10}{c}{$Huber$ $M$} \\ \hline
$\theta$&$0.0$&$  1.0894$&$  0.1644$&$  1.3512$&$12$&$  1.0201$&$  0.0895$&$  1.1299$&$10$\\
$\sigma$&$1.0$&$  1.2912$&$  0.1274$&$  0.2122$&$27$&$  1.3318$&$  0.0893$&$  0.1994$&$19$\\
\hline
&   & & $n=     100$ & & &  &$n=     150$  & &  \\ \hline\hline
& $\tau$ & $\hat{\tau}$ & ${Var}(\hat{\tau})$ & ${MSE}(\hat{\tau})$ & $RE$ & $\hat{\tau}$ & ${Var}(\hat{\tau})$ & ${MSE}(\hat{\tau})$ & $RE$ \\ \hline
\multicolumn{10}{c}{$ESH$} \\ \hline
$\theta$&$0.0$&$ -0.0824$&$  0.0504$&$  0.0572$&$100$&$ -0.0692$&$  0.0458$&$  0.0506$&$100$\\
$\sigma$&$1.0$&$  0.9851$&$  0.0187$&$  0.0189$&$100$&$  0.9636$&$  0.0106$&$  0.0119$&$100$\\
$\varepsilon$&$-0.5$&$ -0.4894$&$  0.0047$&$  0.0049$&$100$&$ -0.4958$&$  0.0028$&$  0.0028$&$100$\\
\hline\hline
\multicolumn{10}{c}{$ESN$} \\ \hline
$\theta$&$0.0$&$  0.3546$&$  0.0472$&$  0.1730$&$33$&$  0.3134$&$  0.0359$&$  0.1341$&$38$\\
$\sigma$&$1.0$&$  1.4866$&$  0.0814$&$  0.3181$&$6$&$  1.5389$&$  0.0716$&$  0.3619$&$3$\\
$\varepsilon$&$-0.5$&$ -0.4603$&$  0.0191$&$  0.0207$&$24$&$ -0.4924$&$  0.0109$&$  0.0110$&$26$\\
\hline\hline
\multicolumn{10}{c}{$ESL$} \\ \hline
$\theta$&$0.0$&$  0.3927$&$  0.0350$&$  0.1892$&$30$&$  0.3832$&$  0.0271$&$  0.1739$&$29$\\
$\sigma$&$1.0$&$  0.6942$&$  0.0059$&$  0.0994$&$19$&$  0.7056$&$  0.0045$&$  0.0912$&$13$\\
$\varepsilon$&$-0.5$&$ -0.3091$&$  0.0109$&$  0.0474$&$10$&$ -0.3166$&$  0.0071$&$  0.0408$&$7$\\
\hline\hline
\multicolumn{10}{c}{$ESt$} \\ \hline
$\theta$&$0.0$&$ -0.0568$&$  0.0594$&$  0.0626$&$91$&$ -0.0796$&$  0.0348$&$  0.0509$&$99$\\
$\sigma$&$1.0$&$  0.7204$&$  0.0046$&$  0.0828$&$23$&$  0.7276$&$  0.0032$&$  0.0774$&$15$\\
$\varepsilon$&$-0.5$&$ -0.5554$&$  0.0218$&$  0.0248$&$20$&$ -0.5638$&$  0.0121$&$  0.0162$&$17$\\
\hline
\multicolumn{10}{c}{$Huber$ $M$} \\ \hline
$\theta$&$0.0$&$  1.0369$&$  0.0646$&$  1.1398$&$5$&$  1.0252$&$  0.0588$&$  1.1099$&$5$\\
$\sigma$&$1.0$&$  1.4452$&$  0.0589$&$  0.2571$&$7$&$  1.2744$&$  0.0780$&$  0.1533$&$8$\\
\hline
\end{tabular}}
\end{table}
\begin{table}[t]
\centering
\label{ESHcomsen0835n}
\caption{Asymmetric M (ESH), ML and M Estimators ($\varepsilon=-0.8$):$c_1=   -0.10, c_2=    6.40$, $k=1.4$}
\scalebox{0.9}{
\begin{tabular}{cccccccccc}\hline
 &   & & $n=      30$ & & &  &$n=      50$  & &  \\ \hline\hline
& $\tau$ & $\hat{\tau}$ & ${Var}(\hat{\tau})$ & ${MSE}(\hat{\tau})$ & $RE$ & $\hat{\tau}$ & ${Var}(\hat{\tau})$ & ${MSE}(\hat{\tau})$ & $RE$ \\ \hline
\multicolumn{10}{c}{$ESH$} \\ \hline
$\theta$&$0.0$&$  0.0689$&$  0.0886$&$  0.0934$&$100$&$ -0.0262$&$  0.0377$&$  0.0383$&$100$\\
$\sigma$&$1.0$&$  1.0331$&$  0.0899$&$  0.0910$&$100$&$  1.0235$&$  0.0578$&$  0.0584$&$100$\\
$\varepsilon$&$-0.8$&$ -0.7178$&$  0.0227$&$  0.0294$&$100$&$ -0.7189$&$  0.0105$&$  0.0171$&$100$\\
\hline\hline
\multicolumn{10}{c}{$ESN$} \\ \hline
$\theta$&$0.0$&$  0.4769$&$  0.1218$&$  0.3492$&$27$&$  0.3285$&$  0.0700$&$  0.1779$&$22$\\
$\sigma$&$1.0$&$  1.7538$&$  0.4441$&$  1.0123$&$9$&$  1.7556$&$  0.2643$&$  0.8352$&$7$\\
$\varepsilon$&$-0.8$&$ -0.6577$&$  0.0311$&$  0.0513$&$57$&$ -0.7121$&$  0.0173$&$  0.0250$&$68$\\
\hline\hline
\multicolumn{10}{c}{$ESL$} \\ \hline
$\theta$&$0.0$&$  0.7885$&$  0.1531$&$  0.7748$&$12$&$  0.7100$&$  0.0721$&$  0.5762$&$7$\\
$\sigma$&$1.0$&$  0.8315$&$  0.0406$&$  0.0690$&$132$&$  0.8008$&$  0.0194$&$  0.0591$&$99$\\
$\varepsilon$&$-0.8$&$ -0.4075$&$  0.0331$&$  0.1872$&$16$&$ -0.4351$&$  0.0181$&$  0.1513$&$11$\\
\hline\hline
\multicolumn{10}{c}{$ESt$} \\ \hline
$\theta$&$0.0$&$  0.1809$&$  0.1671$&$  0.1999$&$47$&$  0.1678$&$  0.0788$&$  0.1069$&$36$\\
$\sigma$&$1.0$&$  0.7402$&$  0.0156$&$  0.0830$&$110$&$  0.7365$&$  0.0101$&$  0.0795$&$73$\\
$\varepsilon$&$-0.8$&$ -0.7586$&$  0.0469$&$  0.0486$&$61$&$ -0.7262$&$  0.0264$&$  0.0319$&$53$\\
\hline
\multicolumn{10}{c}{$Huber$ $M$} \\ \hline
$\theta$&$0.0$&$  1.8830$&$  8.1574$&$  11.7032$&$1$&$  1.8797$&$  6.1185$&$  9.6521$&$0$\\
$\sigma$&$1.0$&$  1.3302$&$  0.8120$&$  0.9209$&$10$&$  1.9495$&$  0.7885$&$  1.6899$&$4$\\
\hline
 &   & & $n=     100$ & & &  &$n=     150$  & &  \\ \hline\hline
& $\tau$ & $\hat{\tau}$ & ${Var}(\hat{\tau})$ & ${MSE}(\hat{\tau})$ & $RE$ & $\hat{\tau}$ & ${Var}(\hat{\tau})$ & ${MSE}(\hat{\tau})$ & $RE$ \\ \hline
\multicolumn{10}{c}{$ESH$} \\ \hline
$\theta$&$0.0$&$ -0.0811$&$  0.0166$&$  0.0232$&$100$&$ -0.0959$&$  0.0060$&$  0.0152$&$100$\\
$\sigma$&$1.0$&$  0.9754$&$  0.0208$&$  0.0214$&$100$&$  0.9858$&$  0.0171$&$  0.0173$&$100$\\
$\varepsilon$&$-0.8$&$ -0.7595$&$  0.0038$&$  0.0065$&$100$&$ -0.7639$&$  0.0028$&$  0.0041$&$100$\\
\hline\hline
\multicolumn{10}{c}{$ESN$} \\ \hline
$\theta$&$0.0$&$  0.2700$&$  0.0371$&$  0.1100$&$21$&$  0.2194$&$  0.0244$&$  0.0726$&$21$\\
$\sigma$&$1.0$&$  1.7370$&$  0.1247$&$  0.6679$&$3$&$  1.8112$&$  0.1056$&$  0.7637$&$2$\\
$\varepsilon$&$-0.8$&$ -0.7397$&$  0.0089$&$  0.0125$&$52$&$ -0.7679$&$  0.0060$&$  0.0070$&$58$\\
\hline\hline
\multicolumn{10}{c}{$ESL$} \\ \hline
$\theta$&$0.0$&$  0.6761$&$  0.0425$&$  0.4996$&$5$&$  0.6810$&$  0.0281$&$  0.4918$&$3$\\
$\sigma$&$1.0$&$  0.7981$&$  0.0086$&$  0.0494$&$43$&$  0.8164$&$  0.0062$&$  0.0399$&$43$\\
$\varepsilon$&$-0.8$&$ -0.4410$&$  0.0097$&$  0.1386$&$5$&$ -0.4415$&$  0.0066$&$  0.1351$&$3$\\
\hline\hline
\multicolumn{10}{c}{$ESt$} \\ \hline
$\theta$&$0.0$&$  0.0500$&$  0.0433$&$  0.0458$&$51$&$  0.0323$&$  0.0223$&$  0.0234$&$65$\\
$\sigma$&$1.0$&$  0.7362$&$  0.0041$&$  0.0737$&$29$&$  0.7427$&$  0.0033$&$  0.0695$&$25$\\
$\varepsilon$&$-0.8$&$ -0.7899$&$  0.0141$&$  0.0142$&$45$&$ -0.8009$&$  0.0076$&$  0.0076$&$54$\\
\hline
\multicolumn{10}{c}{$Huber$ $M$} \\ \hline
$\theta$&$0.0$&$  2.0201$&$  1.8066$&$  5.8874$&$0$&$  2.1522$&$  1.1568$&$  5.7889$&$0$\\
$\sigma$&$1.0$&$  1.8971$&$  0.9565$&$  1.7602$&$1$&$  2.1260$&$  0.6387$&$  1.9067$&$1$\\
\hline
\end{tabular}}
\label{ESHcomsen0835n}
\end{table}

\subsection{Real Data Application on Estimations of Location, Scale and Skewness Parameters}
The real data sets are considerd to show the performance of the asymmetric M-estimation (ESH). As an indicator of the goodness of fitness, $AIC$ (Akaike  information criterion) ve $BIC$ (Bayesian  information criterion) are considered. They are defined as the following forms:

\begin{eqnarray}\label{AICBIC}
AIC&=&2k-2log(L(\hat{\tau};x_1,x_2,...,x_n)) \\
BIC&=&-2log(L(\hat{\tau};x_1,x_2,...,x_n))+klog(n)
\end{eqnarray}

\textbf{Example 1}: The data set in the website http://discover.nci.nih.gov/datasetsNature2000.jsp is analyized by \cite{Arslan09} and \cite{PurHol05}. In this study, the asymmetric M-estimator (ESH), ESN, ESL, ESt, N (Normal) and Huber M-estimation (H) distributions are used to analyze the data set. The tuning constants of asymmetric Huber M-estimators are  $c_1=-0.1,c_2=0.3$, the tuning constant of Huber M-estimators is $k=0.2$ and the parameter $\nu$ that is a tuning constant for the ESt distribution is $\nu=5$. Table \ref{ornekbir} gives the estimates of parameters, log(L), AIC and BIC values.

\begin{table}[t]
\centering
\caption{Example $1$: Estimates of parameters, $logL$, AIC ve BIC}
\label{ornekbir}
\scalebox{0.8}{
\begin{tabular}{c|cccccc}
\hline
              & ESH      & ESN      & ESL      & ESt & N  & H \\ \hline \hline
$\hat{\theta}$      &   0.0386(0.0815)  &  0.1240(0.1677) &   0.2480 &   0.2157(0.1170) &   0.2838(0.0047)  &  0.0332(0.1606)\\
$\hat{\sigma}$      & 0.1195(0.0491) &    0.5139(0.0469)  &  0.3033   &  0.3778(0.0329) &   0.5330(0.0023)  &  0.1355(0.0875)  \\
$\hat{\varepsilon}$ & -0.2049(0.2551)  & -0.1839(0.1884) &  -0.0452 &   -0.0373(0.1788)           &   -      &  - \\
$logL$          & \textbf{ 27.1501} & -15.7852 &   8.2659 & -45.3664 & -17.8779 &  24.1149  \\
AIC           &  \textbf{-48.3002} &  37.5703 & -10.5318  & 96.7328 &   39.7559  & -44.2297 \\
BIC           &  \textbf{-42.0171}  & 43.8534 &   -4.2487 & 103.0159 &   43.9446 & -40.0410\\ \hline
\end{tabular}}
\end{table}

\textbf{Example 2}: The tuning constants of asymmetric Huber M-estimators are  $c_1=-0.25,c_2=0.1$, the tuning constant of Huber M-estimators is $k=0.25$ and the parameter $\nu$ that is a tuning constant for the ESt distribution is $\nu=5$.  Table \ref{ornekiki} gives the estimates of parameters, log(L), AIC and BIC values. The tuning constants were tired until the smallest values of AIC and BIC are gotten for these two examples.

\begin{table}[t]
\centering
\caption{Example $2$: Estimates of parameters, $logL$, AIC ve BIC}
\scalebox{0.8}{
\begin{tabular}{c|cccccc}
\hline
              & ESH      & ESN      & ESL      & ESt & N  & H \\ \hline \hline
$\hat{\theta}$      &   0.1130(0.1287) &  0.2371(0.4414) &  -0.4931 &  -0.1179(0.3899) &  -0.8721(0.0505) &  -0.2550(0.3717)\\
$\hat{\sigma}$      &  0.2260(0.0928)  &  1.5598(0.1423) &   0.9042 &   1.3172(0.4006) &   1.7419(0.0252)   & 0.3837(0.1983) \\
$\hat{\varepsilon}$ &  0.3144(0.2286) & 0.5231(0.1633) &   0.1345 &   0.2950(0.1709)    &   -      &  - \\
$logL$          &   \textbf{-20.5648} & -82.3323&  -57.3653& -113.2674&  -88.9332&  -48.8577  \\
AIC           &   \textbf{47.1297} & 170.6645 & 120.7305 & 232.5347 & 181.8664 & 101.7154 \\
BIC           &    \textbf{53.4127} & 176.9476 & 127.0135 & 238.8178  & 186.0551 &  105.9041 \\ \hline
\end{tabular}}
\label{ornekiki}
\end{table}

\section{Regression Application on Asymmetric $M-$Estimation}
The regression model is considered as an application of asymmetric $M-$estimation.
\begin{flalign}\label{lineerbeta}
~~~~~~~~~~~~~~ y_i = \mathbf{x}_i^T \mathbf{b} + u_i,~~~~ i=1,2,...,n&&
\end{flalign}
\noindent where $y$ is dependent variable. $\mathbf{x}$ is explanatory variable. $b=(b_0,b_1,...,b_{p-1})$ is a vector of parameters. $u$ is error terms. We will get the asymmetric $M-$estimators.

\subsection{Asymmetric $M-$estimation and its Estimators}
The regression model in equation \eqref{lineerbeta} is considered. Suppose that the error terms $u$ are asymmetrically distributed. Then, it can be possible to model the asymmetry in data set via the skewness parameter $\varepsilon$. To estimate the parameters  $\mathbf{b}$, $\sigma$ and $\varepsilon$, the following function $Q$ will be considered.
\begin{equation}\label{minasymobjectivefunctionREG}
 Q( \mathbf{b},\sigma,\varepsilon)= \sum_{i=1}^n   \rho_{ESH} \bigg(\frac{y_i - \mathbf{x}_i^T \mathbf{b}}{\sigma(1-sign(y_i - \mathbf{x}_i^T \mathbf{b})\varepsilon)} \bigg) +nlog(\sigma) \\ \nonumber +\sum_{i=1}^{n} log(1-sign(y_i - \mathbf{x}_i^T \mathbf{b})\varepsilon)
\end{equation}
\noindent The minimization of function $Q$ with respect to parameters will give the estimators of parameters.
\begin{equation}\label{minasymobjectivefunctionREG}
 \frac{\partial Q}{\partial \mathbf{b}}= \sum_{i=1}^n  \psi_{\mathbf{b}} \bigg(\frac{y_i - \mathbf{x}_i^T \mathbf{b}}{\sigma(1-sign(y_i - \mathbf{x}_i^T \mathbf{b})\varepsilon)} \bigg) \frac{\mathbf{x}_i}{\sigma(1-sign(y_i - \mathbf{x}_i^T\mathbf{b})\varepsilon)}=0
\end{equation}

\noindent let $r_i$  be  $\frac{u_i}{\sigma(1-sign(y_i - \mathbf{x}_i^T \mathbf{b})\varepsilon)}$. Then, $w(r_i)=\psi_{\mathbf{b}}(r_i)/r_i$ is defined in the robustness. Then, asymmetric $M-$estimator is

\begin{equation}\label{bestinIF}
 \hat{\mathbf{b}} = [\sum_{i=1}^n \mathbf{x}_i \frac{w_i}{(\hat{\sigma}(1-sign(y_i - \mathbf{x}_i^T \hat{\mathbf{b}})\hat{\varepsilon}))^2} \mathbf{x}_i^T]^{-1}\sum_{i=1}^n \mathbf{x}_i \frac{w_i}{(\hat{\sigma}(1-sign(y_i - \mathbf{x}_i^T \hat{\mathbf{b}})\hat{\varepsilon}))^2}  y_i
\end{equation}
\noindent where $w_i=w(\frac{y_i - \mathbf{x}_i^T \hat{\mathbf{b}}}{\hat{\sigma}(1-sign(y_i - \mathbf{x}_i^T \hat{\mathbf{b}})\hat{\varepsilon})})$ is weight function.

To get the estimator for the parameter $\sigma$, we will get the derivative with respect to $\sigma$
\begin{equation}\label{minasymobjectivefunctionREGsigma}
 \frac{\partial Q}{\partial \sigma} =\frac{-n}{\sigma} + \sum_{i=1}^n  \psi_\sigma \bigg(\frac{y_i - \mathbf{x}_i^T \mathbf{b}}{\sigma(1-sign(y_i - \mathbf{x}_i^T \mathbf{b})\varepsilon)} \bigg) \frac{y_i - \mathbf{x}_i^T \mathbf{b}}{(1-sign(y_i - \mathbf{x}_i^T\mathbf{b})\varepsilon)\sigma^2}=0
\end{equation}
\noindent $r_i =  \frac{u_i}{\sigma(1-sign(y_i - \mathbf{x}_i^T \mathbf{b})\varepsilon)}$. Then, $w(r_i)=\psi_{\sigma}(r_i)/r_i$  is weight function. The asymmetric $M-$estimator of scale parameter is as follow:
\begin{equation}\label{asymregsigma}
\hat{\sigma}^2 = \frac{1}{n} \sum_{i=1}^n w_i  \frac{(y_i - \mathbf{x}_i^T \hat{\mathbf{b}})^2}{(1-sign(y_i - \mathbf{x}_i^T \hat{\mathbf{b}})\hat{\varepsilon})^2}
\end{equation}
\noindent where $w_i=w(\frac{y_i - \mathbf{x}_i^T \hat{\mathbf{b}}}{\hat{\sigma}(1-sign(y_i - \mathbf{x}_i^T \hat{\mathbf{b}})\hat{\varepsilon})})$. The derivative of function $Q( \mathbf{b},\sigma,\varepsilon)$ with respect to parameter $\varepsilon$ is taken
\begin{equation}\label{minasymobjectivefunctionREGepsilon}
 \frac{\partial Q }{\partial \varepsilon} =\frac{sign(y_i - \mathbf{x}_i^T \hat{\mathbf{b}})}{1-sign(y_i - \mathbf{x}_i^T \hat{\mathbf{b}})\varepsilon} - \frac{1}{\sigma} \sum_{i=1}^n  \psi_\varepsilon (r_i) \frac{(y_i - \mathbf{x}_i^T \mathbf{b})sign(y_i - \mathbf{x}_i^T \hat{\mathbf{b}})}{(1-sign(y_i - \mathbf{x}_i^T \hat{\mathbf{b}})\varepsilon)^2}=0
\end{equation}
\noindent  $r_i =  \frac{u_i}{\sigma(1-sign(y_i - \mathbf{x}_i^T \mathbf{b})\varepsilon)}$.  $w(r_i)=\psi_{\varepsilon}(r_i)/r_i$. The asymmetric $M-$estimator of skewness parameter is as follow:
\begin{eqnarray}\label{asymregepsilon}
\hat{\varepsilon}  &=&   \sum_{i=1}^n \bigg[ \frac{sign(y_i - \mathbf{x}_i^T \hat{\mathbf{b}})}{(1-sign(y_i - \mathbf{x}_i^T \hat{\mathbf{b}})\hat{\varepsilon})^2} -  w_i \frac{(y_i - \mathbf{x}_i^T \hat{\mathbf{b}})^2 sign(y_i - \mathbf{x}_i^T \hat{\mathbf{b}})}{\hat{\sigma}^2(1-sign(y_i - \mathbf{x}_i^T \hat{\mathbf{b}})\hat{\varepsilon})^3}\bigg] \nonumber \\
   &/& \sum_{i=1}^n\frac{1}{(1-sign(y_i - \mathbf{x}_i^T \hat{\mathbf{b}})\hat{\varepsilon})^2}
\end{eqnarray}
\noindent $w_i = w(\frac{y_i - \mathbf{x}_i^T \hat{\mathbf{b}}}{\hat{\sigma}(1-sign(y_i - \mathbf{x}_i^T \hat{\mathbf{b}})\hat{\varepsilon})})$.
The weight function for the asymmetric $M-$estimators is as follow
\begin{equation}\label{boundwithcepsilonW}
w(r_i) = \left\{
                \begin{array}{ll}
                  \frac{c_1}{(1+\varepsilon)^2r_i}, & (-\infty,c_1); \\
                  \frac{1}{(1+\varepsilon)^2}, & [c_1,0); \\
                  \frac{1}{(1-\varepsilon)^2}, & [0,c_2]; \\
                  \frac{c_2}{(1-\varepsilon)^2r_i}, & (c_2,\infty).
                \end{array}
              \right.
\end{equation}

\noindent where $r_i=\frac{y_i - \mathbf{x}_i^T \hat{\mathbf{b}}}{\hat{\sigma}(1-sign(y_i - \mathbf{x}_i^T \hat{\mathbf{b}})\hat{\varepsilon})}$. The computation steps are similar to the previous case, because the estimates of parameters location, scale and skewness are considered. Thus, we omitted the steps for the regression case.

\subsection{Simulation Study for Estimators of Regression, Scale and Skewness Parameters}

In the simulation setting, the simulation plan is same with the estimations of location, scale and skewness parameters. The true regression model is as follow:
\begin{equation}\label{lineerbetasimulation}
y_i = 3x_{0i} + 5x_{1i} + x_{2i} - 4x_{3i} + 2x_{4i} - 2 x_{5i}  + u_i,~~~~ i=1,2,...,n
\end{equation}
\noindent The error terms $u$ are distributed asymmetrically. The explanatory variables are $x_0,x_1,x_2,...,x_5$. The inital points of $\mathbf{b}=(b_0,b_1,b_2,b_3,b_4,b_5)$, $\sigma$  and $\varepsilon$ to start the algorithm are the vector of $(0,0,0,0,0,0)$, $MAD$ and $0$, respectively. Three degress of asymmetry are considered. In tables, the results shows that the asymmetric M-estimators outperforms than the maximum likelihood estimators of regression, scale and skewness parameters of ESN, ESL and ESt distributions generally.
\begin{table}[t]
\centering
\caption{Asymmetric M (ESH) and ML Estimators ($\varepsilon=-0.2$): $c_1=   -1.10, c_2=    5.20$}
\label{ESHcomsregen0235n2}
\scalebox{0.9}{
\begin{tabular}{cccccccccc}\hline
 &   & & $n=      30$ & & &  &$n=      50$  & &  \\ \hline\hline
& $\tau$ & $\hat{\tau}$ & ${Var}(\hat{\tau})$ & ${MSE}(\hat{\tau})$ &  $RE$ & $\hat{\tau}$ & ${Var}(\hat{\tau})$ & ${MSE}(\hat{\tau})$ & $RE$ \\ \hline
\multicolumn{10}{c}{$ESH$} \\ \hline
$b_0$&$3.0$&$  2.9675$&$  0.1484$&$  0.1495$&$100$&$  2.8862$&$  0.0841$&$  0.0929$&$100$\\

$b_1$&$5.0$&$  4.9565$&$  0.0344$&$  0.0362$&$100$&$  5.0044$&$  0.0253$&$  0.0253$&$100$\\

$b_2$&$1.0$&$  1.0536$&$  0.0618$&$  0.0647$&$100$&$  0.9990$&$  0.0230$&$  0.0230$&$100$\\

$b_3$&$-4.0$&$ -3.9775$&$  0.0351$&$  0.0356$&$100$&$ -3.9779$&$  0.0260$&$  0.0265$&$100$\\

$b_4$&$2.0$&$  1.9788$&$  0.0458$&$  0.0463$&$100$&$  2.0315$&$  0.0329$&$  0.0339$&$100$\\

$b_5$&$-2.0$&$ -1.9260$&$  0.0457$&$  0.0511$&$100$&$ -1.9911$&$  0.0236$&$  0.0237$&$100$\\

$\sigma$&$1.0$&$  1.0267$&$  0.0135$&$  0.0142$&$100$&$  1.0364$&$  0.0098$&$  0.0112$&$100$\\

$\varepsilon$&$-0.2$&$ -0.1685$&$  0.0141$&$  0.0151$&$100$&$ -0.2003$&$  0.0124$&$  0.0124$&$100$\\
\hline\hline
\multicolumn{10}{c}{$ESN$} \\ \hline

$b_0$&$3.0$&$  3.3478$&$  0.0591$&$  0.1801$&$83$&$  3.3773$&$  0.0219$&$  0.1642$&$57$\\

$b_1$&$5.0$&$  5.0192$&$  0.0807$&$  0.0811$&$45$&$  5.0420$&$  0.0517$&$  0.0535$&$47$\\

$b_2$&$1.0$&$  1.0348$&$  0.0929$&$  0.0941$&$69$&$  1.0235$&$  0.0344$&$  0.0350$&$66$\\

$b_3$&$-4.0$&$ -3.9830$&$  0.0779$&$  0.0781$&$46$&$ -3.9705$&$  0.0473$&$  0.0482$&$55$\\

$b_4$&$2.0$&$  1.9189$&$  0.1256$&$  0.1322$&$35$&$  2.0516$&$  0.0344$&$  0.0371$&$91$\\

$b_5$&$-2.0$&$ -1.9313$&$  0.0697$&$  0.0745$&$69$&$ -1.9883$&$  0.0339$&$  0.0340$&$70$\\

$\sigma$&$1.0$&$  1.1339$&$  0.0783$&$  0.0963$&$15$&$  1.2276$&$  0.1005$&$  0.1523$&$7$\\

$\varepsilon$&$-0.2$&$ -0.0354$&$  0.0068$&$  0.0339$&$44$&$ -0.0810$&$  0.0089$&$  0.0231$&$54$\\
\hline\hline
\multicolumn{10}{c}{$ESL$} \\ \hline

$b_0$&$3.0$&$  3.1078$&$  0.1612$&$  0.1728$&$87$&$  3.1346$&$  0.0803$&$  0.0998$&$93$\\

$b_1$&$5.0$&$  4.7403$&$  0.1761$&$  0.2435$&$15$&$  4.6403$&$  0.1550$&$  0.2844$&$8$\\

$b_2$&$1.0$&$  0.8516$&$  0.0902$&$  0.1122$&$58$&$  0.9095$&$  0.0978$&$  0.1060$&$22$\\

$b_3$&$-4.0$&$ -3.7824$&$  0.2043$&$  0.2517$&$14$&$ -3.7205$&$  0.1581$&$  0.2363$&$11$\\

$b_4$&$2.0$&$  1.7326$&$  0.1880$&$  0.2595$&$18$&$  1.8225$&$  0.0721$&$  0.1071$&$32$\\

$b_5$&$-2.0$&$ -1.8504$&$  0.1347$&$  0.1571$&$33$&$ -1.8451$&$  0.0905$&$  0.1145$&$21$\\

$\sigma$&$1.0$&$  0.6916$&$  0.0384$&$  0.1335$&$11$&$  0.7090$&$  0.0209$&$  0.1056$&$11$\\

$\varepsilon$&$-0.2$&$ -0.1064$&$  0.0194$&$  0.0282$&$53$&$ -0.1351$&$  0.0122$&$  0.0164$&$76$\\

\hline\hline
\multicolumn{10}{c}{$ESt$} \\ \hline
$b_0$&$3.0$&$  3.1056$&$  0.3187$&$  0.3298$&$45$&$  2.9574$&$  0.2092$&$  0.2111$&$44$\\

$b_1$&$5.0$&$  4.9588$&$  0.0862$&$  0.0879$&$41$&$  5.0119$&$  0.0388$&$  0.0390$&$65$\\

$b_2$&$1.0$&$  1.0232$&$  0.1097$&$  0.1103$&$59$&$  1.0038$&$  0.0334$&$  0.0334$&$69$\\

$b_3$&$-4.0$&$ -4.0213$&$  0.0663$&$  0.0668$&$53$&$ -3.9777$&$  0.0398$&$  0.0403$&$66$\\

$b_4$&$2.0$&$  1.9617$&$  0.0981$&$  0.0996$&$46$&$  1.9986$&$  0.0354$&$  0.0354$&$96$\\

$b_5$&$-2.0$&$ -1.9408$&$  0.0581$&$  0.0616$&$83$&$ -1.9755$&$  0.0394$&$  0.0400$&$59$\\

$\sigma$&$1.0$&$  0.7318$&$  0.2703$&$  0.3422$&$4$&$  0.6817$&$  0.0086$&$  0.1100$&$10$\\

$\varepsilon$&$-0.2$&$ -0.1841$&$  0.1730$&$  0.1733$&$9$&$ -0.2764$&$  0.1022$&$  0.1080$&$12$\\

\hline
\end{tabular}}
\end{table}

\begin{table}[t]
\centering
\caption{Asymmetric M (ESH) and ML Estimators  ($\varepsilon=-0.2$): $c_1=   -1.10, c_2=    5.20$}
\scalebox{0.9}{
\begin{tabular}{cccccccccc}\hline
 &   & & $n=     100$ & & &  &$n=     150$  & &  \\ \hline\hline
& $\tau$ & $\hat{\tau}$ & ${Var}(\hat{\tau})$ & ${MSE}(\hat{\tau})$ &  $RE$ & $\hat{\tau}$ & ${Var}(\hat{\tau})$ & ${MSE}(\hat{\tau})$ & $RE$ \\ \hline
\multicolumn{10}{c}{$ESH$} \\ \hline
$b_0$&$3.0$&$  2.7998$&$  0.0428$&$  0.0828$&$100$&$  2.8049$&$  0.0311$&$  0.0692$&$100$\\

$b_1$&$5.0$&$  5.0042$&$  0.0099$&$  0.0099$&$100$&$  4.9814$&$  0.0094$&$  0.0098$&$100$\\

$b_2$&$1.0$&$  0.9832$&$  0.0066$&$  0.0069$&$100$&$  0.9984$&$  0.0068$&$  0.0068$&$100$\\

$b_3$&$-4.0$&$ -4.0094$&$  0.0144$&$  0.0144$&$100$&$ -3.9802$&$  0.0065$&$  0.0069$&$100$\\

$b_4$&$2.0$&$  2.0029$&$  0.0081$&$  0.0081$&$100$&$  2.0063$&$  0.0064$&$  0.0064$&$100$\\

$b_5$&$-2.0$&$ -1.9821$&$  0.0092$&$  0.0095$&$100$&$ -1.9935$&$  0.0079$&$  0.0079$&$100$\\

$\sigma$&$1.0$&$  1.0473$&$  0.0044$&$  0.0066$&$100$&$  1.0447$&$  0.0025$&$  0.0045$&$100$\\

$\varepsilon$&$-0.2$&$ -0.2070$&$  0.0044$&$  0.0044$&$100$&$ -0.2157$&$  0.0039$&$  0.0041$&$100$\\
\hline\hline
\multicolumn{10}{c}{$ESN$} \\ \hline

$b_0$&$3.0$&$  3.3573$&$  0.0224$&$  0.1501$&$55$&$  3.3669$&$  0.0093$&$  0.1439$&$48$\\

$b_1$&$5.0$&$  5.0270$&$  0.0187$&$  0.0194$&$51$&$  4.9870$&$  0.0181$&$  0.0183$&$53$\\

$b_2$&$1.0$&$  0.9920$&$  0.0137$&$  0.0138$&$50$&$  1.0008$&$  0.0096$&$  0.0096$&$71$\\

$b_3$&$-4.0$&$ -4.0387$&$  0.0267$&$  0.0282$&$51$&$ -3.9871$&$  0.0094$&$  0.0095$&$72$\\

$b_4$&$2.0$&$  2.0043$&$  0.0215$&$  0.0215$&$38$&$  2.0051$&$  0.0114$&$  0.0114$&$56$\\

$b_5$&$-2.0$&$ -1.9846$&$  0.0155$&$  0.0157$&$60$&$ -1.9977$&$  0.0131$&$  0.0131$&$60$\\

$\sigma$&$1.0$&$  1.2661$&$  0.0389$&$  0.1097$&$6$&$  1.2943$&$  0.0254$&$  0.1120$&$4$\\

$\varepsilon$&$-0.2$&$ -0.0951$&$  0.0051$&$  0.0161$&$27$&$ -0.0980$&$  0.0046$&$  0.0150$&$27$\\
\hline\hline
\multicolumn{10}{c}{$ESL$} \\ \hline

$b_0$&$3.0$&$  3.0491$&$  0.0971$&$  0.0996$&$83$&$  3.0047$&$  0.0839$&$  0.0840$&$82$\\

$b_1$&$5.0$&$  4.6588$&$  0.1551$&$  0.2815$&$3$&$  4.6069$&$  0.1201$&$  0.2746$&$4$\\

$b_2$&$1.0$&$  0.8519$&$  0.0355$&$  0.0679$&$10$&$  0.9336$&$  0.0409$&$  0.0453$&$15$\\

$b_3$&$-4.0$&$ -3.6718$&$  0.0844$&$  0.1921$&$8$&$ -3.6635$&$  0.0846$&$  0.1908$&$4$\\

$b_4$&$2.0$&$  1.8236$&$  0.0774$&$  0.1065$&$8$&$  1.8640$&$  0.0835$&$  0.1020$&$6$\\

$b_5$&$-2.0$&$ -1.8170$&$  0.0808$&$  0.1143$&$8$&$ -1.8302$&$  0.0596$&$  0.0884$&$9$\\

$\sigma$&$1.0$&$  0.7443$&$  0.0227$&$  0.0881$&$8$&$  0.7855$&$  0.0212$&$  0.0673$&$7$\\

$\varepsilon$&$-0.2$&$ -0.1175$&$  0.0079$&$  0.0147$&$30$&$ -0.1789$&$  0.0130$&$  0.0135$&$31$\\

\hline\hline
\multicolumn{10}{c}{$ESt$} \\ \hline
$b_0$&$3.0$&$  2.9470$&$  0.0871$&$  0.0982$&$84$&$  2.8777$&$  0.0823$&$  0.0973$&$71$\\

$b_1$&$5.0$&$  4.9989$&$  0.0134$&$  0.0134$&$74$&$  4.9699$&$  0.0128$&$  0.0130$&$75$\\

$b_2$&$1.0$&$  0.9859$&$  0.0109$&$  0.0111$&$62$&$  1.0080$&$  0.0092$&$  0.0093$&$73$\\

$b_3$&$-4.0$&$ -4.0174$&$  0.0184$&$  0.0187$&$77$&$ -3.9771$&$  0.0076$&$  0.0081$&$85$\\

$b_4$&$2.0$&$  2.0119$&$  0.0162$&$  0.0163$&$50$&$  2.0060$&$  0.0068$&$  0.0069$&$93$\\

$b_5$&$-2.0$&$ -1.9831$&$  0.0142$&$  0.0145$&$66$&$ -2.0023$&$  0.0127$&$  0.0127$&$63$\\

$\sigma$&$1.0$&$  0.7386$&$  0.0042$&$  0.0726$&$9$&$  0.7685$&$  0.0030$&$  0.0566$&$8$\\

$\varepsilon$&$-0.2$&$ -0.2604$&$  0.0497$&$  0.0513$&$9$&$ -0.2935$&$  0.0332$&$  0.0419$&$10$\\

\hline
\end{tabular}}
\end{table}
\begin{table}[t]
\centering
\caption{DAsymmetric M (ESH) and ML Estimators ($\varepsilon=-0.5$): $c_1=   -0.30, c_2=    5.30$}
\scalebox{0.9}{
\begin{tabular}{cccccccccc}\hline
 &   & & $n=      30$ & & &  &$n=      50$  & &  \\ \hline\hline
& $\tau$ & $\hat{\tau}$ & ${Var}(\hat{\tau})$ & ${MSE}(\hat{\tau})$ &  $RE$ & $\hat{\tau}$ & ${Var}(\hat{\tau})$ & ${MSE}(\hat{\tau})$ & $RE$ \\ \hline
\multicolumn{10}{c}{$ESH$} \\ \hline
$b_0$&$3.0$&$  3.0629$&$  0.1043$&$  0.1083$&$100$&$  2.9531$&$  0.0506$&$  0.0528$&$100$\\

$b_1$&$5.0$&$  5.0174$&$  0.0496$&$  0.0499$&$100$&$  5.0219$&$  0.0193$&$  0.0197$&$100$\\

$b_2$&$1.0$&$  1.0027$&$  0.0419$&$  0.0419$&$100$&$  1.0106$&$  0.0163$&$  0.0164$&$100$\\

$b_3$&$-4.0$&$ -4.0016$&$  0.0511$&$  0.0511$&$100$&$ -4.0097$&$  0.0266$&$  0.0267$&$100$\\

$b_4$&$2.0$&$  2.0080$&$  0.0472$&$  0.0473$&$100$&$  2.0126$&$  0.0122$&$  0.0123$&$100$\\

$b_5$&$-2.0$&$ -1.9778$&$  0.0446$&$  0.0451$&$100$&$ -2.0076$&$  0.0171$&$  0.0172$&$100$\\

$\sigma$&$1.0$&$  1.0301$&$  0.1086$&$  0.1095$&$100$&$  0.9721$&$  0.0100$&$  0.0108$&$100$\\

$\varepsilon$&$-0.5$&$ -0.5097$&$  0.0278$&$  0.0279$&$100$&$ -0.5230$&$  0.0089$&$  0.0094$&$100$\\
\hline\hline
\multicolumn{10}{c}{$ESN$} \\ \hline

$b_0$&$3.0$&$  3.9324$&$  0.0643$&$  0.9338$&$12$&$  3.9961$&$  0.0563$&$  1.0485$&$5$\\

$b_1$&$5.0$&$  5.0071$&$  0.0898$&$  0.0899$&$56$&$  5.0391$&$  0.0630$&$  0.0645$&$31$\\

$b_2$&$1.0$&$  1.0402$&$  0.0853$&$  0.0869$&$48$&$  1.0118$&$  0.0603$&$  0.0605$&$27$\\

$b_3$&$-4.0$&$ -4.0389$&$  0.0957$&$  0.0972$&$53$&$ -4.0706$&$  0.0783$&$  0.0833$&$32$\\

$b_4$&$2.0$&$  2.0174$&$  0.0877$&$  0.0880$&$54$&$  2.0291$&$  0.0631$&$  0.0639$&$19$\\

$b_5$&$-2.0$&$ -2.0281$&$  0.1258$&$  0.1266$&$36$&$ -2.0227$&$  0.0709$&$  0.0714$&$24$\\

$\sigma$&$1.0$&$  1.4340$&$  0.3739$&$  0.5623$&$19$&$  1.5196$&$  0.2371$&$  0.5071$&$2$\\

$\varepsilon$&$-0.5$&$ -0.1376$&$  0.0126$&$  0.1439$&$19$&$ -0.1717$&$  0.0098$&$  0.1176$&$8$\\
\hline\hline
\multicolumn{10}{c}{$ESL$} \\ \hline

$b_0$&$3.0$&$  4.3295$&$  0.1289$&$  0.5785$&$19$&$  3.6835$&$  0.1095$&$  0.5767$&$9$\\

$b_1$&$5.0$&$  4.7333$&$  0.1948$&$  0.2659$&$19$&$  4.7566$&$  0.1906$&$  0.2499$&$8$\\

$b_2$&$1.0$&$  0.9364$&$  0.1432$&$  0.1473$&$28$&$  0.8366$&$  0.0845$&$  0.1112$&$15$\\

$b_3$&$-4.0$&$ -3.7703$&$  0.2556$&$  0.3084$&$17$&$ -3.7676$&$  0.1708$&$  0.2248$&$12$\\

$b_4$&$2.0$&$  1.9062$&$  0.1776$&$  0.1864$&$25$&$  1.9135$&$  0.1227$&$  0.1302$&$9$\\

$b_5$&$-2.0$&$ -1.9351$&$  0.1649$&$  0.1691$&$27$&$ -1.8726$&$  0.1472$&$  0.1634$&$11$\\

$\sigma$&$1.0$&$  0.7930$&$  0.0520$&$  0.0949$&$115$&$  0.7978$&$  0.0513$&$  0.0922$&$12$\\

$\varepsilon$&$-0.5$&$ -0.1640$&$  0.0205$&$  0.1334$&$21$&$ -0.1658$&$  0.0114$&$  0.1231$&$8$\\

\hline\hline
\multicolumn{10}{c}{$ESt$} \\ \hline
$b_0$&$3.0$&$  3.0856$&$  0.1987$&$  0.2060$&$53$&$  3.0557$&$  0.1434$&$  0.1465$&$36$\\

$b_1$&$5.0$&$  4.9587$&$  0.0886$&$  0.0903$&$55$&$  5.0196$&$  0.0340$&$  0.0344$&$57$\\

$b_2$&$1.0$&$  0.9923$&$  0.0832$&$  0.0832$&$50$&$  1.0168$&$  0.0378$&$  0.0381$&$43$\\

$b_3$&$-4.0$&$ -4.0140$&$  0.1452$&$  0.1454$&$35$&$ -4.0211$&$  0.0310$&$  0.0315$&$85$\\

$b_4$&$2.0$&$  1.9574$&$  0.1279$&$  0.1297$&$36$&$  2.0144$&$  0.0359$&$  0.0361$&$34$\\

$b_5$&$-2.0$&$ -2.0362$&$  0.0874$&$  0.0887$&$51$&$ -1.9977$&$  0.0250$&$  0.0250$&$69$\\

$\sigma$&$1.0$&$  0.6531$&$  0.0376$&$  0.1580$&$69$&$  0.6898$&$  0.0104$&$  0.1066$&$10$\\

$\varepsilon$&$-0.5$&$ -0.5929$&$  0.0841$&$  0.0927$&$30$&$ -0.5810$&$  0.0680$&$  0.0746$&$13$\\

\hline
\end{tabular}}
\end{table}
\begin{table}[t]
\centering
\caption{Asymmetric M (ESH) and ML Estimators ($\varepsilon=-0.5$): $c_1=   -0.30, c_2=    5.30$}
\scalebox{0.9}{
\begin{tabular}{cccccccccc}\hline
 &   & & $n=     100$ & & &  &$n=     150$  & &  \\ \hline\hline
& $\tau$ & $\hat{\tau}$ & ${Var}(\hat{\tau})$ & ${MSE}(\hat{\tau})$ &  $RE$ & $\hat{\tau}$ & ${Var}(\hat{\tau})$ & ${MSE}(\hat{\tau})$ & $RE$ \\ \hline
\multicolumn{10}{c}{$ESH$} \\ \hline
$b_0$&$3.0$&$  2.8901$&$  0.0155$&$  0.0285$&$100$&$  2.8595$&$  0.0083$&$  0.0280$&$100$\\

$b_1$&$5.0$&$  5.0040$&$  0.0039$&$  0.0039$&$100$&$  4.9915$&$  0.0033$&$  0.0034$&$100$\\

$b_2$&$1.0$&$  1.0027$&$  0.0053$&$  0.0053$&$100$&$  0.9958$&$  0.0041$&$  0.0041$&$100$\\

$b_3$&$-4.0$&$ -3.9976$&$  0.0037$&$  0.0037$&$100$&$ -3.9993$&$  0.0029$&$  0.0029$&$100$\\

$b_4$&$2.0$&$  1.9966$&$  0.0025$&$  0.0025$&$100$&$  2.0003$&$  0.0014$&$  0.0014$&$100$\\

$b_5$&$-2.0$&$ -2.0032$&$  0.0036$&$  0.0036$&$100$&$ -2.0024$&$  0.0024$&$  0.0024$&$100$\\

$\sigma$&$1.0$&$  0.9736$&$  0.0038$&$  0.0045$&$100$&$  0.9815$&$  0.0028$&$  0.0032$&$100$\\

$\varepsilon$&$-0.5$&$ -0.5164$&$  0.0027$&$  0.0030$&$100$&$ -0.5072$&$  0.0016$&$  0.0016$&$100$\\
\hline\hline
\multicolumn{10}{c}{$ESN$} \\ \hline

$b_0$&$3.0$&$  3.9450$&$  0.0240$&$  0.9171$&$3$&$  3.9202$&$  0.0151$&$  0.8619$&$3$\\

$b_1$&$5.0$&$  5.0383$&$  0.0316$&$  0.0331$&$12$&$  5.0192$&$  0.0206$&$  0.0210$&$16$\\

$b_2$&$1.0$&$  0.9864$&$  0.0287$&$  0.0289$&$18$&$  1.0054$&$  0.0171$&$  0.0171$&$24$\\

$b_3$&$-4.0$&$ -4.0314$&$  0.0328$&$  0.0338$&$11$&$ -4.0368$&$  0.0171$&$  0.0184$&$16$\\

$b_4$&$2.0$&$  1.9962$&$  0.0234$&$  0.0234$&$11$&$  2.0151$&$  0.0148$&$  0.0150$&$9$\\

$b_5$&$-2.0$&$ -2.0215$&$  0.0220$&$  0.0225$&$16$&$ -2.0180$&$  0.0107$&$  0.0111$&$22$\\

$\sigma$&$1.0$&$  1.5933$&$  0.1231$&$  0.4751$&$1$&$  1.6149$&$  0.0927$&$  0.4807$&$1$\\

$\varepsilon$&$-0.5$&$ -0.2053$&$  0.0074$&$  0.0942$&$3$&$ -0.2249$&$  0.0066$&$  0.0823$&$2$\\
\hline\hline
\multicolumn{10}{c}{$ESL$} \\ \hline

$b_0$&$3.0$&$  3.5221$&$  0.1011$&$  0.3737$&$8$&$  3.4696$&$  0.0808$&$  0.3014$&$9$\\

$b_1$&$5.0$&$  4.5857$&$  0.3390$&$  0.5106$&$1$&$  4.5264$&$  0.1515$&$  0.3757$&$1$\\

$b_2$&$1.0$&$  0.8813$&$  0.0979$&$  0.1120$&$5$&$  0.9225$&$  0.0603$&$  0.0663$&$6$\\

$b_3$&$-4.0$&$ -3.6655$&$  0.1793$&$  0.2912$&$1$&$ -3.6455$&$  0.1292$&$  0.2549$&$1$\\

$b_4$&$2.0$&$  1.8401$&$  0.0709$&$  0.0965$&$3$&$  1.8302$&$  0.0605$&$  0.0893$&$2$\\

$b_5$&$-2.0$&$ -1.8937$&$  0.1010$&$  0.1119$&$3$&$ -1.8087$&$  0.0618$&$  0.1009$&$2$\\

$\sigma$&$1.0$&$  0.8658$&$  0.1095$&$  0.1275$&$4$&$  0.8737$&$  0.0444$&$  0.0603$&$5$\\

$\varepsilon$&$-0.5$&$ -0.2029$&$  0.0118$&$  0.1000$&$3$&$ -0.2169$&$  0.0107$&$  0.0909$&$2$\\

\hline\hline
\multicolumn{10}{c}{$ESt$} \\ \hline
$b_0$&$3.0$&$  2.9438$&$  0.0577$&$  0.0608$&$47$&$  2.9242$&$  0.0342$&$  0.0400$&$70$\\

$b_1$&$5.0$&$  4.9897$&$  0.0192$&$  0.0192$&$20$&$  4.9879$&$  0.0099$&$  0.0100$&$34$\\

$b_2$&$1.0$&$  0.9879$&$  0.0157$&$  0.0158$&$33$&$  0.9934$&$  0.0075$&$  0.0075$&$55$\\

$b_3$&$-4.0$&$ -3.9983$&$  0.0127$&$  0.0127$&$30$&$ -4.0147$&$  0.0091$&$  0.0093$&$31$\\

$b_4$&$2.0$&$  1.9915$&$  0.0091$&$  0.0092$&$27$&$  2.0051$&$  0.0053$&$  0.0053$&$26$\\

$b_5$&$-2.0$&$ -2.0110$&$  0.0085$&$  0.0086$&$41$&$ -1.9994$&$  0.0074$&$  0.0074$&$32$\\

$\sigma$&$1.0$&$  0.7626$&$  0.0052$&$  0.0615$&$7$&$  0.7793$&$  0.0039$&$  0.0526$&$6$\\

$\varepsilon$&$-0.5$&$ -0.5857$&$  0.0205$&$  0.0279$&$11$&$ -0.5775$&$  0.0113$&$  0.0173$&$9$\\

\hline
\end{tabular}}
\end{table}
\begin{table}[t]
\centering
\caption{Asymmetric M (ESH) and ML Estimators ($\varepsilon=-0.8$): $c_1=   -0.01, c_2=    6.20$}
\scalebox{0.9}{
\begin{tabular}{cccccccccc}\hline
 &   & & $n=      30$ & & &  &$n=      50$  & &  \\ \hline\hline
& $\tau$ & $\hat{\tau}$ & ${Var}(\hat{\tau})$ & ${MSE}(\hat{\tau})$ &  $RE$ & $\hat{\tau}$ & ${Var}(\hat{\tau})$ & ${MSE}(\hat{\tau})$ & $RE$ \\ \hline
\multicolumn{10}{c}{$ESH$} \\ \hline
$b_0$&$3.0$&$  3.4040$&$  0.2139$&$  0.3771$&$100$&$  3.1880$&$  0.1209$&$  0.1562$&$100$\\

$b_1$&$5.0$&$  4.9738$&$  0.1024$&$  0.1031$&$100$&$  5.0228$&$  0.0768$&$  0.0773$&$100$\\

$b_2$&$1.0$&$  1.0093$&$  0.1496$&$  0.1496$&$100$&$  1.0040$&$  0.0468$&$  0.0468$&$100$\\

$b_3$&$-4.0$&$ -4.0286$&$  0.1756$&$  0.1764$&$100$&$ -3.9875$&$  0.0475$&$  0.0476$&$100$\\

$b_4$&$2.0$&$  1.8981$&$  0.1859$&$  0.1963$&$100$&$  2.0424$&$  0.0417$&$  0.0435$&$100$\\

$b_5$&$-2.0$&$ -1.9979$&$  0.1245$&$  0.1245$&$100$&$ -1.9866$&$  0.0307$&$  0.0308$&$100$\\

$\sigma$&$1.0$&$  1.1666$&$  0.3282$&$  0.3559$&$100$&$  1.0377$&$  0.1754$&$  0.1768$&$100$\\

$\varepsilon$&$-0.8$&$ -0.7395$&$  0.0463$&$  0.0500$&$100$&$ -0.8319$&$  0.0223$&$  0.0233$&$100$\\
\hline\hline
\multicolumn{10}{c}{$ESN$} \\ \hline

$b_0$&$3.0$&$  4.6000$&$  0.2060$&$ 2.7660$&$14$&$  4.5671$&$  0.0794$&$  2.5350$&$6$\\

$b_1$&$5.0$&$  5.0640$&$  0.2286$&$  0.2327$&$44$&$  5.0713$&$  0.1417$&$  0.1468$&$53$\\

$b_2$&$1.0$&$  1.0032$&$  0.1768$&$  0.1769$&$85$&$  1.0373$&$  0.1224$&$  0.1238$&$38$\\

$b_3$&$-4.0$&$ -4.0691$&$  0.2639$&$  0.2687$&$66$&$ -4.0339$&$  0.0845$&$  0.0857$&$56$\\

$b_4$&$2.0$&$  1.9688$&$  0.2828$&$  0.2837$&$69$&$  2.0409$&$  0.1182$&$  0.1199$&$36$\\

$b_5$&$-2.0$&$ -2.0679$&$  0.2640$&$  0.2686$&$46$&$ -2.0494$&$  0.0867$&$  0.0892$&$35$\\

$\sigma$&$1.0$&$  2.0024$&$  0.8743$&$  1.8791$&$19$&$  2.1712$&$  0.7601$&$  2.1317$&$8$\\

$\varepsilon$&$-0.8$&$ -0.1903$&$  0.0095$&$  0.3812$&$13$&$ -0.2500$&$  0.0109$&$  0.3134$&$7$\\
\hline\hline
\multicolumn{10}{c}{$ESL$} \\ \hline

$b_0$&$3.0$&$  4.1359$&$  0.2345$&$  1.5248$&$25$&$  4.0596$&$  0.1143$&$  1.2372$&$13$\\

$b_1$&$5.0$&$  4.7404$&$  0.3929$&$  0.4602$&$22$&$  4.7188$&$  0.1857$&$  0.2647$&$29$\\

$b_2$&$1.0$&$  1.0128$&$  0.1411$&$  0.1412$&$106$&$  0.9760$&$  0.1342$&$  0.1348$&$35$\\

$b_3$&$-4.0$&$ -3.7965$&$  0.3274$&$  0.3688$&$48$&$ -3.6758$&$  0.1485$&$  0.2536$&$19$\\

$b_4$&$2.0$&$  1.8841$&$  0.2825$&$  0.2960$&$66$&$  1.8651$&$  0.1080$&$  0.1262$&$35$\\

$b_5$&$-2.0$&$ -1.9028$&$  0.2800$&$  0.2894$&$43$&$ -1.9026$&$  0.1127$&$  0.1222$&$25$\\

$\sigma$&$1.0$&$  0.9353$&$  0.1219$&$  0.1260$&$282$&$  0.9547$&$  0.0549$&$  0.0570$&$310$\\

$\varepsilon$&$-0.8$&$ -0.1981$&$  0.0150$&$  0.3773$&$13$&$ -0.2408$&$  0.0118$&$  0.3245$&$7$\\

\hline\hline
\multicolumn{10}{c}{$ESt$} \\ \hline
$b_0$&$3.0$&$  3.2988$&$  0.5966$&$  0.6859$&$55$&$  3.2092$&$  0.1182$&$  0.1619$&$96$\\

$b_1$&$5.0$&$  4.8167$&$  0.5956$&$  0.6292$&$17$&$  4.9793$&$  0.0493$&$  0.0497$&$155$\\

$b_2$&$1.0$&$  0.9902$&$  0.4074$&$  0.4075$&$37$&$  0.9896$&$  0.0538$&$  0.0539$&$87$\\

$b_3$&$-4.0$&$ -3.7916$&$  0.6577$&$  0.7012$&$25$&$ -3.9680$&$  0.0785$&$  0.0795$&$60$\\

$b_4$&$2.0$&$  1.8722$&$  0.2551$&$  0.2714$&$72$&$  2.0209$&$  0.0895$&$  0.0899$&$48$\\

$b_5$&$-2.0$&$ -1.8194$&$  0.4842$&$  0.5168$&$24$&$ -2.0189$&$  0.0477$&$  0.0480$&$64$\\

$\sigma$&$1.0$&$  0.6930$&$  0.0436$&$  0.1379$&$258$&$  0.7431$&$  0.0278$&$  0.0938$&$189$\\

$\varepsilon$&$-0.8$&$ -0.6726$&$  0.0726$&$  0.0889$&$56$&$ -0.7778$&$  0.0458$&$  0.0463$&$50$\\

\hline
\end{tabular}}
\end{table}
\begin{table}[t]
\centering
\caption{Asymmetric M (ESH) and ML Estimators ($\varepsilon=-0.8$): $c_1=   -0.01, c_2=    6.20$}
\label{ESHcomsregen0815n2}
\scalebox{0.9}{
\begin{tabular}{cccccccccc}\hline
 &   & & $n=     100$ & & &  &$n=     150$  & &  \\ \hline\hline
& $\tau$ & $\hat{\tau}$ & ${Var}(\hat{\tau})$ & ${MSE}(\hat{\tau})$ &  $RE$ & $\hat{\tau}$ & ${Var}(\hat{\tau})$ & ${MSE}(\hat{\tau})$ & $RE$ \\ \hline
\multicolumn{10}{c}{$ESH$} \\ \hline

$b_0$&$3.0$&$  3.0883$&$  0.0105$&$  0.0183$&$100$&$  2.9272$&$  0.0092$&$  0.0145$&$100$\\

$b_1$&$5.0$&$  5.0099$&$  0.0081$&$  0.0082$&$100$&$  5.0020$&$  0.0012$&$  0.0012$&$100$\\

$b_2$&$1.0$&$  0.9911$&$  0.0115$&$  0.0116$&$100$&$  1.0028$&$  0.0013$&$  0.0014$&$100$\\

$b_3$&$-4.0$&$ -4.0027$&$  0.0060$&$  0.0061$&$100$&$ -4.0009$&$  0.0022$&$  0.0022$&$100$\\

$b_4$&$2.0$&$  1.9883$&$  0.0063$&$  0.0064$&$100$&$  1.9970$&$  0.0023$&$  0.0023$&$100$\\

$b_5$&$-2.0$&$ -1.9932$&$  0.0047$&$  0.0047$&$100$&$ -1.9987$&$  0.0016$&$  0.0016$&$100$\\

$\sigma$&$1.0$&$  0.9697$&$  0.0077$&$  0.0086$&$100$&$  0.9668$&$  0.0034$&$  0.0045$&$100$\\

$\varepsilon$&$-0.8$&$ -0.7452$&$  0.0042$&$  0.0072$&$100$&$ -0.7639$&$  0.0045$&$  0.0058$&$100$\\
\hline\hline
\multicolumn{10}{c}{$ESN$} \\ \hline

$b_0$&$3.0$&$  3.5584$&$  0.3346$&$  0.5631$&$3$&$  3.5236$&$  0.2240$&$  0.4760$&$3$\\

$b_1$&$5.0$&$  5.0601$&$  0.0545$&$  0.0581$&$14$&$  5.0768$&$  0.0218$&$  0.0277$&$4$\\

$b_2$&$1.0$&$  1.0179$&$  0.0397$&$  0.0400$&$29$&$  1.0039$&$  0.0320$&$  0.0321$&$4$\\

$b_3$&$-4.0$&$ -4.0586$&$  0.0418$&$  0.0452$&$13$&$ -4.0589$&$  0.0317$&$  0.0352$&$6$\\

$b_4$&$2.0$&$  2.0107$&$  0.0409$&$  0.0410$&$16$&$  2.0366$&$  0.0191$&$  0.0205$&$11$\\

$b_5$&$-2.0$&$ -2.0415$&$  0.0511$&$  0.0529$&$9$&$ -2.0498$&$  0.0295$&$  0.0320$&$5$\\

$\sigma$&$1.0$&$  2.1763$&$  0.3331$&$  1.7168$&$1$&$  2.1919$&$  0.2047$&$  1.6253$&$0$\\

$\varepsilon$&$-0.8$&$ -0.2917$&$  0.0068$&$  0.2652$&$2$&$ -0.3090$&$  0.0052$&$  0.2462$&$2$\\
\hline\hline
\multicolumn{10}{c}{$ESL$} \\ \hline

$b_0$&$3.0$&$  3.9831$&$  0.1062$&$  1.0727$&$2$&$  3.9723$&$  0.1267$&$  1.0721$&$1$\\

$b_1$&$5.0$&$  4.6342$&$  0.1360$&$  0.2698$&$3$&$  4.6800$&$  0.1799$&$  0.2563$&$1$\\

$b_2$&$1.0$&$  0.9216$&$  0.0817$&$  0.0878$&$13$&$  0.9271$&$  0.0729$&$  0.0782$&$2$\\

$b_3$&$-4.0$&$ -3.6813$&$  0.1384$&$  0.2400$&$3$&$ -3.6560$&$  0.1091$&$  0.2274$&$1$\\

$b_4$&$2.0$&$  1.8291$&$  0.0965$&$  0.1257$&$5$&$  1.8276$&$  0.0551$&$  0.0848$&$3$\\

$b_5$&$-2.0$&$ -1.8124$&$  0.0778$&$  0.1130$&$4$&$ -1.8239$&$  0.0724$&$  0.1035$&$2$\\

$\sigma$&$1.0$&$  0.9879$&$  0.0429$&$  0.0430$&$20$&$  0.9977$&$  0.0327$&$  0.0327$&$14$\\

$\varepsilon$&$-0.8$&$ -0.2682$&$  0.0125$&$  0.2953$&$2$&$ -0.2683$&$  0.0116$&$  0.2942$&$2$\\

\hline\hline
\multicolumn{10}{c}{$ESt$} \\ \hline
$b_0$&$3.0$&$  3.0865$&$  0.0394$&$  0.0469$&$39$&$  3.0720$&$  0.0270$&$  0.0322$&$45$\\

$b_1$&$5.0$&$  5.0086$&$  0.0081$&$  0.0082$&$100$&$  5.0022$&$  0.0046$&$  0.0046$&$27$\\

$b_2$&$1.0$&$  0.9993$&$  0.0108$&$  0.0108$&$107$&$  0.9977$&$  0.0031$&$  0.0031$&$43$\\

$b_3$&$-4.0$&$ -4.0001$&$  0.0063$&$  0.0063$&$97$&$ -3.9974$&$  0.0054$&$  0.0054$&$40$\\

$b_4$&$2.0$&$  1.9980$&$  0.0082$&$  0.0082$&$78$&$  1.9895$&$  0.0058$&$  0.0059$&$39$\\

$b_5$&$-2.0$&$ -1.9989$&$  0.0073$&$  0.0073$&$64$&$ -2.0069$&$  0.0046$&$  0.0046$&$35$\\

$\sigma$&$1.0$&$  0.7816$&$  0.0055$&$  0.0532$&$16$&$  0.7993$&$  0.0037$&$  0.0440$&$10$\\

$\varepsilon$&$-0.8$&$ -0.8105$&$  0.0146$&$  0.0147$&$49$&$ -0.7969$&$  0.0082$&$  0.0082$&$71$\\

\hline
\end{tabular}}
\end{table}

\subsection{Real Data Application on Estimations of Regression, Scale and Skewness Parameters}
\textbf{Example 3}: The Martin–Marietta data set was analyezed by \cite{ButlerMc90}, \cite{AzCap03,DicMon04,Genc07,Arslan09,ArslanGenc09} and \cite{Acitasetal13}. The estimates of the regression parameters were obitaned by these studies. They assume that th The data set consists of 60 monthly observations from January 1982 to December 1986. \cite{ButlerMc90} introduce a linear regression model $y=b_0 + b_1CRSP + u$ where $y$ is the excess rate of the Martin–Marietta company, $x =$ CRSP is an index of the excess rate of return for the New York market, and  $u$ is an error term. The tuning constants of asymmetric Huber M-estimators are  $c_1=-0.015,c_2=0.03$, the tuning constant of Huber M-estimators is $k=0.03$ and the parameter $\nu$ that is a tuning constant for the ESt distribution is $\nu=1.5$. Table \ref{table:realexample1linesreg} gives the estimates of parameters, log(L), AIC and BIC values.

\begin{table}[t]
\centering
\caption{Example $3$: Estimates of parameters, $logL$, AIC ve BIC}
\label{table:realexample1linesreg}
\scalebox{0.9}{
\begin{tabular}{c|cccccc}
\hline
                    & ESH      & ESN      & ESL      & ESt            & N  & H \\ \hline \hline
$\hat{b}_0$         &    -0.0047 &  -0.0177 &  -0.0099&   -0.0092   & 0.0011 &  -0.0016 \\
$\hat{b}_1$         &     0.3607  &  1.5118 &   0.6846 &   1.0009  &  1.8025  &  0.4241 \\
$\hat{\sigma}$      &     0.1240  &  0.1105 &   0.0546 &   0.0633  &  0.1210  &  0.1373 \\
$\hat{\varepsilon}$ &   -0.1116   &-0.2581  & -0.1093  & -0.0459   &      -   &      - \\
$logL$          &    \textbf{70.1209}  & 58.8884 &  63.7627 &   66.1233 &  53.9302  & 63.0189 \\
AIC           & \textbf{-134.2418} & -111.7768 & -121.5255 & -126.2466 & -101.8604 & -120.0377 \\
BIC           &  \textbf{-127.9588} & -105.4938 & -115.2424 & -119.9636 & -95.5774 & -113.7547 \\ \hline
\end{tabular}}
\end{table}

\textbf{Example 4}:  The same regression model considered in the previous example $3$ is again taken. This data set can get the PET package. It is called as "la". The dependent variable "y=la\$PET510" and the explanatory variable  "x=la\$bflow" are standardized. "bflow" is a variable representing the blood flow. The variable "PET510" is measured via the PET (Positron Emission Tomography) machines. The detailed discussion can be found in the package named as gamlss.nl. The \cite{Langeetal89,JonFad03,RigStasin06} studies also analyzed this data set. The tuning constants of asymmetric Huber M-estimators are  $c_1=-1,c_2=0.7$, the tuning constant of Huber M-estimators is $k=1$ and the parameter $\nu$ that is a tuning constant for the ESt distribution is $\nu=2$. If a value is bigger than the $Q_3 + 1.5 IQR=1.6686$, the added value with $y$ direction is considered to be an outlier. The maximum value of explanatory variable is $5.4276$, the added value is $5$. Here, $Q_3$ is the third quantile. $IQR$ is an interquartile range. The sample size is $n=251$. After adding one outlier, the sample size is $n=252$.  Table \ref{table:realexample2linesregPET} and \ref{table:realexample2linesregPETOUT} give the estimates of parameters, log(L), AIC and BIC values.

\begin{table}[t]
\centering
\caption{Example $3$: Estimates of parameters, $logL$, AIC ve BIC}
\label{table:realexample2linesregPET}
\scalebox{0.9}{
\begin{tabular}{c|cccccc}
\hline
                    & ESH      & ESN      & ESL      & ESt            & N  & H \\ \hline \hline
$\hat{b}_0$         &  -0.0773 &  -0.0841 &  -0.0280 &   0.0118 &  -0.0000 &  -0.0229 \\
$\hat{b}_1$         &    0.7299 &   0.8094  &  0.4586 &   0.6881 &   0.8046 &   0.7688 \\
$\hat{\sigma}$      &   0.6013  &  0.9618 &   0.4908 &   0.5486 &   0.9980 &    0.6629  \\
$\hat{\varepsilon}$ &  -0.0320  & -0.1599  & 0.0000  &  0.0813    &    -   &      -  \\
$logL$          &    \textbf{-176.1575} & -265.8107 & -258.1030& -223.6345 &-274.4067 &-191.9880  \\
AIC           &    \textbf{360.3149} & 539.6215 &524.2060&  455.2690&  554.8134 & 389.9761 \\
BIC           &      \textbf{374.4167}&  553.7233 &538.3078 & 469.3708  &565.3897&  400.5524 \\ \hline
\end{tabular}}
\end{table}

\begin{table}[t]
\centering
\caption{Example $4$(Added outlier): Estimates of parameters, $logL$, AIC ve BIC}
\label{table:realexample2linesregPETOUT}
\scalebox{0.9}{
\begin{tabular}{c|cccccc}
\hline
                    & ESH      & ESN      & ESL      & ESt            & N  & H \\ \hline \hline
$\hat{b}_0$         &    -0.0876  & -0.0850 &  -0.0316   & 0.0103   & 0.0288 &  -0.0175 \\
$\hat{b}_1$         &           0.7440 &   0.9262 &   0.4746  & 0.6930   & 0.9491  &  0.7959\\
$\hat{\sigma}$      &       0.6305 &   1.1451   & 0.5220  &  0.5538&    1.2504 &   0.6943 \\
$\hat{\varepsilon}$ &    -0.0528  & -0.2532 &  -0.0323  &  0.0734  &                  -   &      - \\
$logL$          &    \textbf{-189.6646}& -314.7171& -273.7550 &-231.9491 &-334.5094 &-207.4125 \\
AIC           &    \textbf{387.3292}  &637.4342 & 555.5100 & 471.8983 & 675.0188  &420.8250 \\
BIC           &    \textbf{401.4469}  &651.5519 & 569.6277 & 486.0160&  685.6070  &431.4133 \\ \hline
\multicolumn{7}{l}{}
\end{tabular}
}
\end{table}

\clearpage
\section{Conclusions}\label{conclusions}

Asymmetric Huber M-estimators are suggested by using the $-log(f)$ as an objective function in robustness. Asymmetric M-estimators have a skewness parameter to model the potential skewness in data set. The asymptotic properties of asymmetric M-estimators are examined. Firstly, the existence and uniqueness of the proposed objective function with respect to parameters are examined. After that, the asymptotic normality of estimators can be shown via the well known Taylor expansion of the proposed asymmetric M-function and the multivariate Slutsky's lemma is used. For the proposed estimators, the asymptotic variance-covariance matrix is provided. The influence function as a local robustness property of estimators was provided. The breakdown point as an indicator of global robustness of location estimator is shown to be $1/2$ after providing the assumptions of the robustness criteria of breakdown. An application on regression was also considered. Real data examples for both cases were provided. The results show the competence of our proposed estimators when there is a potential asymmetry in data set. The asymptotic properties of regression case  will be considered as a comprehensive study. The asymmetric M-estimates in R software will be added. The asymmetric forms of Welsch, Hampel functions will be suggested to model the asymmetry in data set.


\end{document}